\documentclass[
  11pt     %
]{article}
\usepackage{amssymb}
\usepackage{amsfonts}
\usepackage{}

\usepackage{color}
\definecolor{black}{rgb}{0,0,0}

\definecolor{red}{rgb}{1,0,0}

\definecolor{blue}{rgb}{0,0,1}

\usepackage{amstext,amsthm,amssymb,amsmath}
\usepackage{bm}
\usepackage{graphicx}
\usepackage{subfigure}
\usepackage{wrapfig}
\usepackage{fullpage}
\usepackage{color}
\usepackage{multirow}
\usepackage{tabulary}
\usepackage{booktabs}
\usepackage{enumerate}
\usepackage{epstopdf}

\newcommand{\average}[1]{\{\!\!\!\{{#1}\}\!\!\!\}}
\newcommand{\jump}[1]{[\![{#1}]\!]}

\newtheorem{theorem}{Theorem}[section]

\newtheorem{remark}[theorem]{Remark}

\newcommand{\norm}[1]{\left\|#1\right\|}

\hfuzz=\maxdimen
\title{Sparse Generalized Multiscale Finite Element Methods and their applications}

\author{Eric Chung\thanks{Department of Mathematics, The Chinese University of Hong Kong, Shatin, Hong Kong SAR. Eric Chung's research is partially supported
by CUHK Direct Grant for Research 2014-15.}
\and Yalchin Efendiev\thanks{Department of Mathematics  and Institute for Scientific Computation (ISC),
Texas A\&M University,
College Station, Texas 77843-3368, USA} \and Wing Tat Leung\thanks{Department of Mathematics  and Institute for Scientific Computation (ISC),
Texas A\&M University,
College Station, Texas 77843-3368, USA. This work is partially supported by the U.S. Department of Energy Office of Science, Office of Advanced Scientific Computing Research, Applied Mathematics program under Award Number DE-FG02-13ER26165} \and Guanglian Li\thanks{Institute for Numerical Simulation, the University of Bonn,
Wegelerstrasse 6, 53115 Bonn, Germany}}

\begin{document}

\maketitle

\begin{abstract}

In a number of previous papers \cite{egh12, eglmsMSDG, eglp13oversampling,calo2014multiscale,chung2014adaptive1,chung2015generalizedperforated,chung2015residual,chung2015online}, local (coarse grid)
multiscale
model reduction techniques are developed
using a Generalized Multiscale Finite Element
Method. In these approaches,
multiscale basis functions are constructed using local snapshot spaces,
where a snapshot space is a large space that represents the solution
behavior in a coarse block. In a number of applications (e.g., those discussed in the paper),
one may have a sparsity in the snapshot space for an appropriate choice
of a snapshot space. More precisely,
the solution may only involve a portion of the snapshot space.
In this case, one can use sparsity techniques (\cite{candes2006compressive,candes2008restricted,candes2006robust,mackey2014compressive,schaeffer2013sparse, jiao2014primal}) to identify multiscale basis
functions. In this paper, we consider two such sparse local
multiscale model reduction approaches.

In the first approach (which is used for parameter-dependent multiscale PDEs),
we use local minimization techniques, such as sparse POD, to identify multiscale
basis functions, which are sparse in the snapshot space.
These minimization techniques use $l_1$ minimization to find
local multiscale basis functions, which are further used for finding
the solution. In the second approach (which is used
for the Helmholtz equation), we directly apply $l_1$
minimization techniques to solve the underlying PDEs. This approach is more
expensive as it involves a large snapshot space; however, in this example,
we can not identify a local minimization principle, such as local generalized
SVD.

All our numerical results assume the sparsity and we discuss this assumption
for the snapshot spaces.
Moreover, we discuss the computational savings provided by our approach.
The sparse solution allows a fast evaluation of stiffness matrices
and downscaling the solution to the fine grid since
the reduced dimensional solution representation is sparse
in terms of local snapshot vectors.
 Numerical results are presented, which show the convergence
of the proposed method and the sparsity of the solution.

\end{abstract}

\section{Introduction}

\subsection{Multiscale problems and the problem of sparsity}

Simulations of multiscale problems are expensive and,  typically, require
some type of a model reduction.
Our approaches seek adaptive reduced-order models, locally in space,
and construct multiscale basis functions in each coarse region
to represent the solution space. These approaches share
common concepts with homogenization and upscaling methods
\cite{papanicolau1978asymptotic, bakhvalov1989homogenisation, eh09, weh02,numerical-homo,2d-waves, fish2008mathematical, fish2013practical},
where local effective properties are constructed.
In contrast, in multiscale methods, local multiscale
basis functions \cite{hw97,jennylt03,melenk1996partition, eh09, oz07,GMsFEM-wave,GMsFEM-mixed,elastic-jcp,GMsFEM-elastic,Chung-Leung-Cicp, fish2013practical, zch15, hl15}
are constructed to represent the solution space.
These basis functions are typically constructed in the snapshot spaces
\cite{egh12}.
In this paper, we investigate cases, when the basis functions are sparse
in the snapshot space.
We discuss several examples, which include multiscale
parameter-dependent problems and the Helmholtz equation.
The parameter-dependent multiscale problems are motivated by stochastic
problems, where the parameter is used to describe the uncertainties.

\subsection{Sparse GMsFEM Concepts}

In this paper, we use the GMsFEM framework
(\cite{galvis2015generalized,Ensemble, eglmsMSDG, eglp13oversampling,calo2014multiscale,chung2014adaptive1,chung2015generalizedperforated,chung2015residual,chung2015online,egh12})
 and investigate
the sparsity within GMsFEM snapshots. To illustrate the main idea
of our approach, we consider
\[
L u = f,
\]
where $L$ is a differential operator. For example, in the paper, we consider
parameter-dependent heterogeneous flows, $Lu=-div(\kappa(x;\mu)\nabla u)$,
 and the Helmholtz equation, $Lu=-div(\kappa(x)\nabla u) - \Omega^2 n(x) u$.
The main idea of GMsFEM
is to construct a snapshot space and identify a subspace,
called the offline or online space depending whether the problem is
parameter-dependent.
This subspace is used
 to solve the underlying problem at a reduced cost.
The snapshot and online spaces are constructed in each coarse element
(see next section for more precise definitions), where a coarse element is a region,
which is much larger than the characteristic fine-length scale
(see Figure \ref{snapshotoverview1}).
 For each coarse region,
$\tau_j$, we construct snapshot vectors, $\{ \psi^j_i \}$
(here $i$ is the numbering of the snapshot functions),
that represent the local solution space.
We denote the snapshot space
by
\[
V_{\text{snap}}^{\tau_j}=\text{Span}_i\{\psi_i^j \}, \ \
V_{\text{snap}}=\text{Span}_{i,j}\{\psi_i^j \}.
\]
In GMsFEM, the online spaces are constructed using the elements of local
snapshot functions.
In many examples, the snapshot space can be large and
the online space can be a sparse subspace of the snapshot space.
The objective of this paper is to investigate these cases.

\subsection{Snapshot spaces}



The snapshot spaces play an important role in the GMsFEM. They are designed
to capture the solution space locally and are used to preserve
some features of the solution space, e.g., mass conservation.
Typical snapshot spaces consist of local solutions constructed using
some sets of boundary conditions or right hand sides.
With an appropriate choice of snapshot spaces (e.g., using oversampling
\cite{eglp13oversampling}), one can improve the convergence
of GMsFEM substantially.

\begin{figure}[htb]
  \centering
   \includegraphics[width=3in, height=3in]{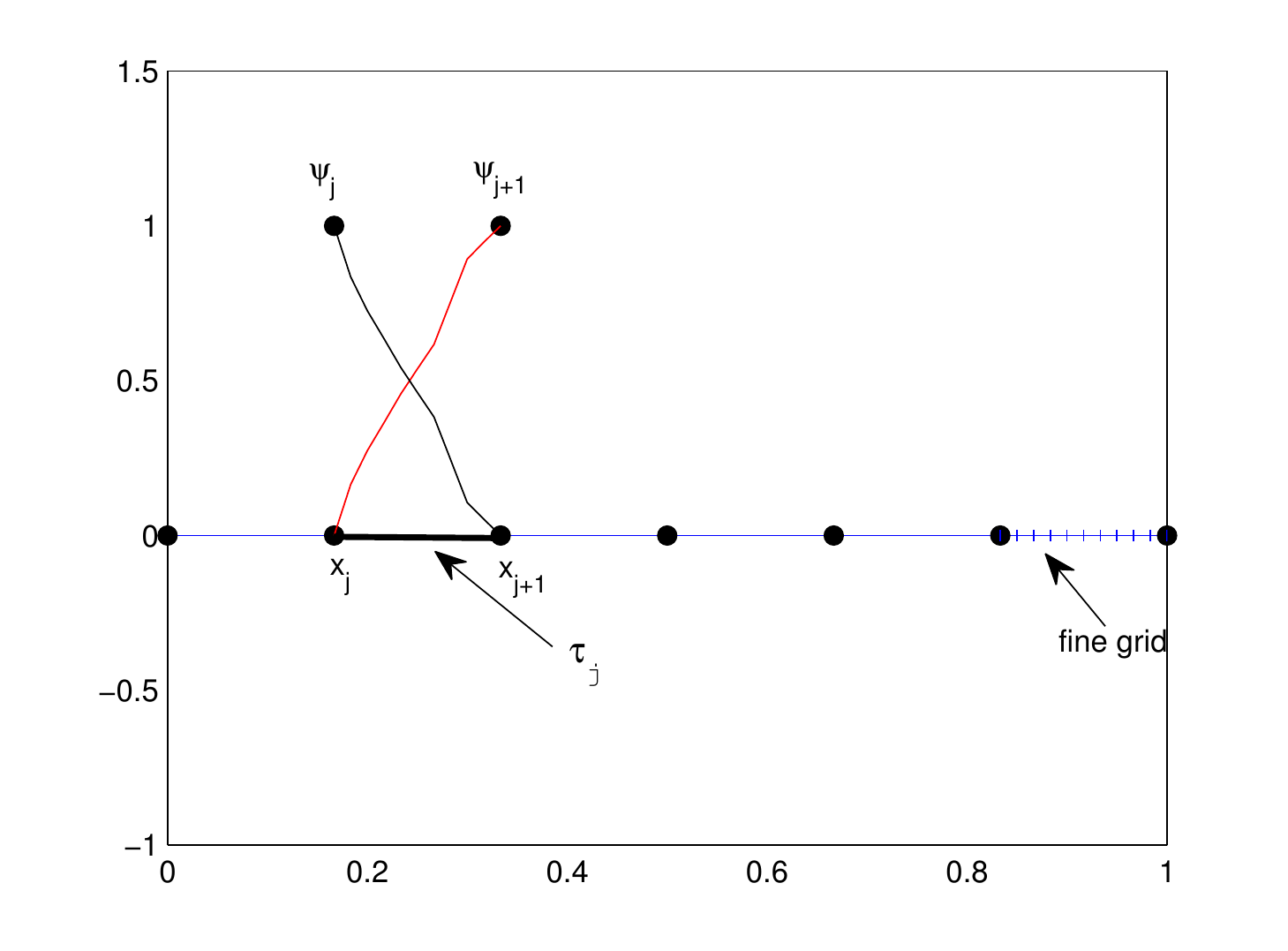}
  \caption{Illustration of snapshot concepts in one dimensional example.}
  \label{snapshot1D}
\end{figure}
\begin{figure}[htb]
  \centering
  \includegraphics[width=3in, height=2in]{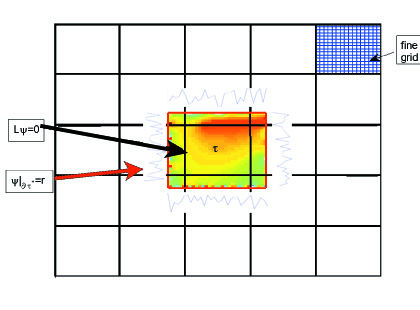}
  \includegraphics[width=3in, height=2in]{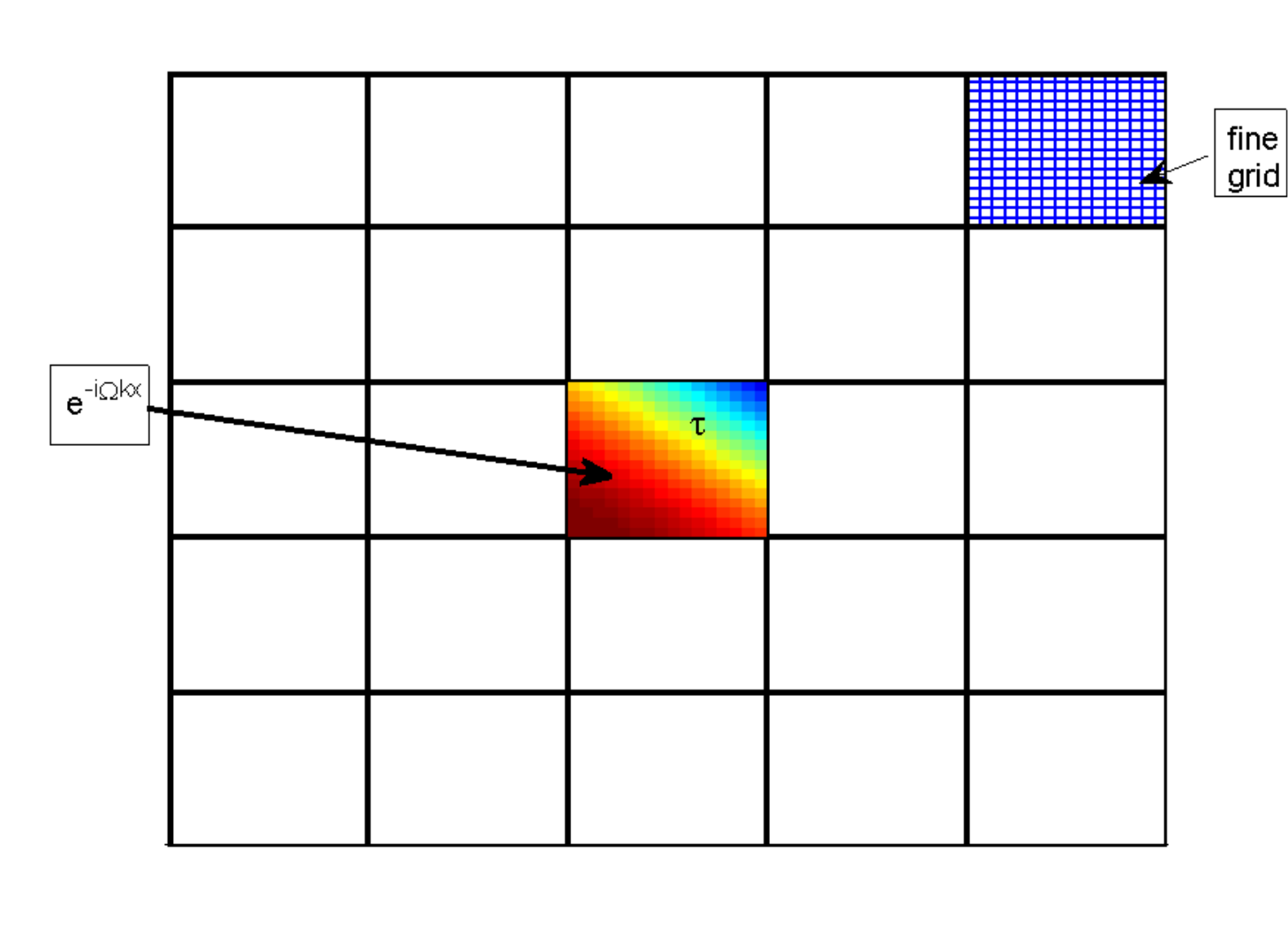}
  \caption{Illustration of snapshot concepts. Left: For the first problem; Right: For the second problem.}
  \label{snapshotoverview1}
\end{figure}

To convey the concept of snapshot spaces, we present some examples.
We will consider two examples discussed above. We start with a simplified
 example related to the parameter-dependent case,
 $L_0u=-div(\kappa(x;\mu=0)\nabla u)$, i.e., the problem without a parameter.
In each coarse-grid block
$\tau_j$ (see the left plot in Figure \ref{snapshotoverview1}),
we consider a local solution
\begin{equation}
\label{def:snap1}
L_0(\psi_i^j)=0\ \text{in}\ \tau_j
\end{equation}
subject to some boundary conditions, where these boundary conditions
play an important role in defining snapshot functions.
One option is to choose all possible boundary
conditions considering all unit vectors on the boundary of $\tau_j$. More precisely,
$\psi_i^j(x)=\delta_i(x)$ on $\partial \tau_j$, where
$\delta_i(x)$ is $1$ at the node
$i$  and zero elsewhere. The computations of these snapshot functions
 are
expensive. Instead,
we use the boundary conditions, which are randomly distributed numbers
on the fine-grid
nodes of the boundary $\partial \tau_j$
 (see the left plot in Figure \ref{snapshotoverview1}).
The random boundary conditions allow extracting the essential information
provided we choose several more snapshot vectors than the number of modes,
we would like to use.
For {\it parameter-dependent problems}, the snapshot vectors are defined as above
(\ref{def:snap1}) for each  pre-selected value of $\mu_m$.
For example, for {\it one-dimensional case} $\tau_j=[x_j,x_{j+1}]$,
for parameter-independent problem,
the snapshot space in $\tau_j$
consists of two solutions $\psi^j$ and $\psi^{j+1}$, such that
$\psi^n(x_l)=\delta_{nl}$, $n=j,j+1$, $l=j,j+1$,
where $\delta_{jl}$ is the Kronecker symbol and $\psi^n$ ($n=j,j+1$)
 is a solution of
${d\over dx}\left( \kappa(x;\mu=0) {d\over dx}\psi^n\right)=0$ in $\tau_j$
(see Figure \ref{snapshot1D} for illustration).
For {\it parameter-dependent problems},
the snapshot vectors in $\tau_j=[x_j,x_{j+1}]$ are the solutions of
\[
{d\over dx}\left( \kappa(x;\mu_m) {d\over dx}\psi_m^n\right)=0\ \text{in} \ \tau_j,
\]
 $\psi^n_m(x_l)=\delta_{nl}$, $n=j,j+1$, $l=j,j+1$.
For multi-dimensional examples, we can construct the snapshots similarly for
each $\tau_j$ and for using different boundary conditions and
 use one index to represent the snapshot
vectors as $\psi_i^j$.
For the second example,
$Lu=-div(\kappa(x)\nabla u) - \Omega^2 n(x) u$, we choose the snapshot
vectors to be functions $ e^{i \Omega k_i\cdot x}$ for a set of pre-defined values
of $k_i$ on a unit circle (see the right plot in Figure \ref{snapshotoverview1}).

The multiscale basis functions are constructed in the snapshot space. Our earlier
approaches seek a small dimensional subspace of the snapshot space by performing
a local spectral decomposition (based on analysis). However, these approaches
use all snapshot vectors when seeking multiscale basis functions. In a number of
applications, the solution is sparse in the snapshot space. I.e., in the expansion
\[
u=\sum_{i,j} c_{i,j} \psi_i^j,
\]
many coefficients $c_{i,j}$ are zeros. In this case, one can save computational effort
by employing sparsity techniques.  In this paper,
our main goal is to discuss how GMsFEM can be designed if the solution
is sparse in the snapshot space. We describe two classes of approaches
and present a framework for constructing sparse GMsFEM.

The main challenge in these applications is to construct a snapshot space,
where the solution is sparse.
In our first example, this can be achieved,
because an online parameter value $\mu$ can be close to some of the pre-selected
 offline
values of $\mu$'s, and thus, the multiscale basis functions (and the solution)
can have a sparse representation in the snapshot space. In our second example,
we select cases where the solution $u$ contains only a few snapshot
vectors corresponding to directions $k_i$.
We note that if the snapshot space is not chosen carefully, one may not have
the sparsity.
In general, there can be many other examples
and our goal is to show how local multiscale model reduction techniques can
be used for such problems.

\subsection{Approaches for identifying sparse solutions in snapshot spaces}

\begin{figure}[htb]
  \centering
  \includegraphics[width=0.65 \textwidth]{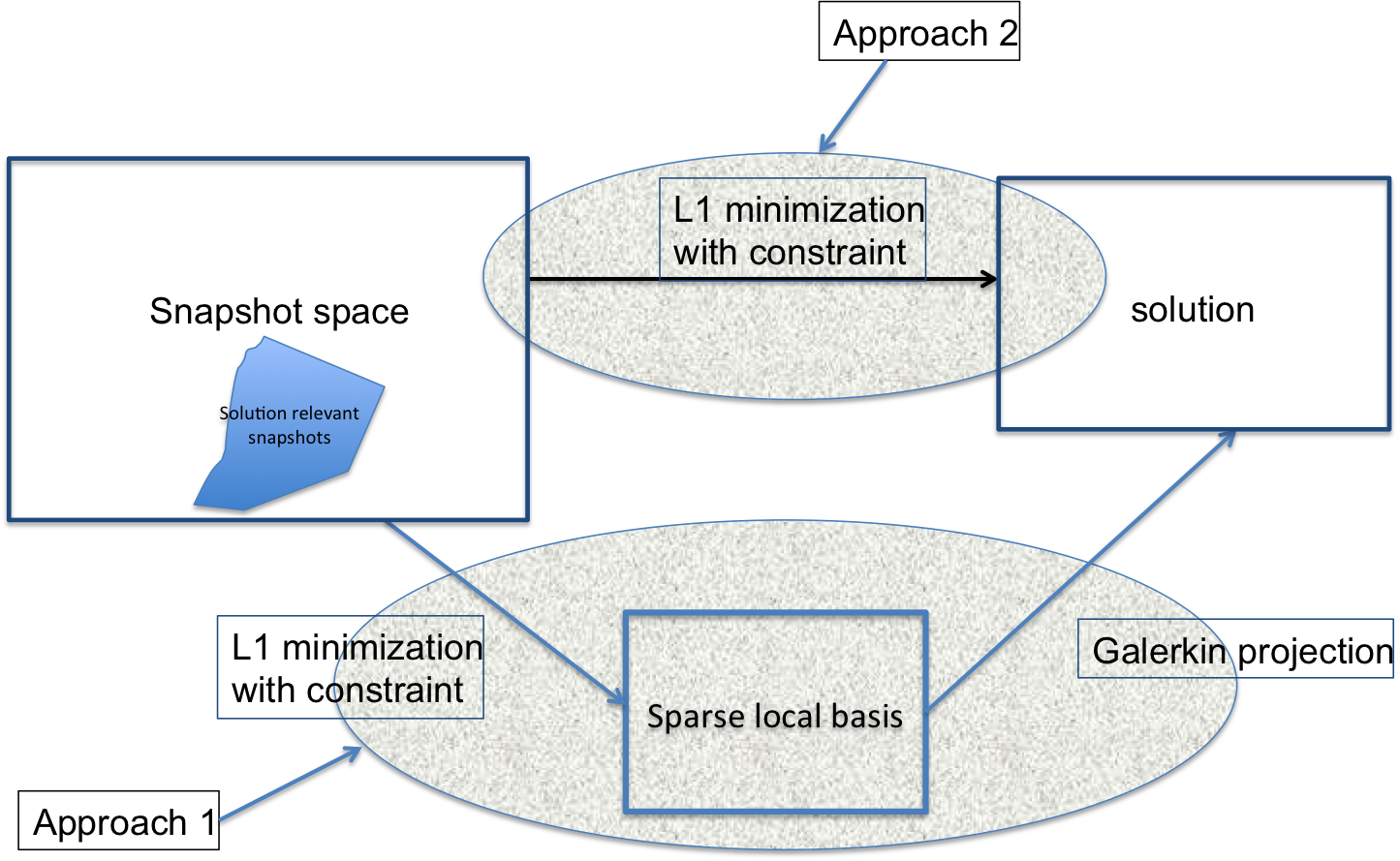}
  \caption{Illustration of our approaches.}
  \label{overview1}
\end{figure}

The main goal
of this paper is to discuss how to explore sparsity ideas within GMsFEM
for constructing
local multiscale basis functions. We consider two distinct cases.
\begin{itemize}

\item First approach: ``Local-Sparse Snapshot Subspace Approach''. Determining the online sparse space locally via local spectral sparse decomposition in the snapshot space (motivated by
parameter-dependent problems).

\item  Second approach: ``Sparse Snapshot Subspace Approach''.  Determining the online space globally via a global solve  (motivated by using plane wave snapshot vectors and the Helmholtz equation).

\end{itemize}
See Figure \ref{overview1} for illustration. We use sparsity techniques
(e.g., \cite{candes2006compressive,candes2008restricted,candes2006robust,mackey2014compressive,schaeffer2013sparse})
to identify local multiscale basis functions
and solve the global problem.

In above approaches, the snapshot functions can be linearly dependent.
In fact, in general, we would like to have a large snapshot space that can
contain
a sparse representation of the solution. In both approaches formulated above,
the linear
dependency is removed. In the first approach, it is removed by sparse POD.
We note that in the original GMsFEM approach \cite{egh12},
POD across all
snapshot functions is used to remove
linear dependence. However, this can result in a loss of sparsity, i.e.,
the solution may contain many nonzero coefficients when represented
in the snapshot space. Thus, for sparsity, it is important to avoid
a POD step across all snapshot vectors. In our first example,
we need to avoid using POD for all $\mu$'s. For this reason, we
design a special sparse POD method using randomized snapshot functions.
It both eliminates linearly dependent snapshot functions and
identifies a sparse solution space. In the second method, $l_1$
minimization can
be used for all snapshots, even if they are linearly dependent.
By adding more snapshot vectors, we hope to identify a sparse representation
of the solution.
The second method will eliminate the linear dependence and identify a sparse
solution. Note that in our example, the snapshot vectors are linearly
independent.

For the first case,
we consider parameter-dependent elliptic equations of the form
\begin{equation} \label{eq:original_NL}
-\mbox{div} \big( \kappa(x;\mu) \, \nabla u  \big)=f \, \, \text{in} \, D,
\end{equation}
where $u=g$ on $\partial D$, and $\mu$ is a parameter.
Some of the existing approaches for parameter-dependent problems closely related to the proposed approaches include
reduced basis techniques \cite{ barrault, ct_POD_DEIM11}. In these approaches, the
reduced order model is constructed via a greedy algorithm.
In the proposed approach, we attempt to approximate
the solution space locally in each coarse block using $l_1$ minimization.
Local multiscale basis functions are constructed using GMsFEM.
In previous approaches, we attempted to compress
the local solutions corresponding to some pre-selected values, $\mu$, in the offline stage. This can lead to large dimensional
offline spaces.  In this paper, we propose an approach to compute eigenvectors using
$l_1$ minimization based on randomized snapshots and oversampling.
This provides an efficient approach to identify sparse eigenvector
representation in GMsFEM.
The proposed approach gives a sparse representation of
the multiscale basis functions in terms of the snapshot space vectors
and has several advantages. (1) It allows quickly assembling of the stiffness
matrix in the online space since it involves a few elements of the
snapshot space. (2) We can downscale the solution to the fine grid
much faster using sparse representation. (3) It avoids a POD based
step proposed in the original GMsFEM formulation (see \cite{egh12}),
which is performed across all $\mu$'s and
which can result to a large dimensional representation
of the online multiscale basis functions in terms of snapshot functions.

In the second example, we consider the Helmholtz equation
\begin{equation} \label{eq:Helmholtz}
-\mbox{div} \big( \kappa(x) \, \nabla u  \big) - \Omega^2 n(x) u =f \, \, \text{in} \, D,
\end{equation}
where $\Omega$ is the frequency.
We will use plane waves as the
snapshot vectors.
In the computational examples considered in this paper, the solution has a few dominant propagating directions,
and is, therefore, spanned by only a few plane waves. This observation leads to the solution sparsity in our snapshot space.
However, the choice of local spectral decomposition
is not available for determining these dominant directions and we study using
$l_1$ minimization directly in the space of snapshot vectors.
We consider not very large snapshot spaces, and snapshot vectors having closed form formulas for constant media properties ($\kappa(x)$ and $n(x)$).
For the Helmholtz equation (\ref{eq:Helmholtz})
with low frequencies, we can expect a sparsity in our examples, which we exploit.
Thus, sparsity techniques are the natural methodologies in these situations.
In general, we can also consider the frequency to be a parameter,
and apply similar techniques used for the first example.


\subsection{Summary of numerical results}

Two test cases for the first approach are presented,
 where we consider parametrized
conductivity fields. In the first case, the conductivity
is parametrized as an affine combination of two heterogeneous
conductivity fields. In the second case, we use a nonlinear
parameter dependence. In particular, we consider an initial conductivity
field with a channel and inclusions. The parametrization is introduced
such that these high-conductivity features spatially move within the domain.
This is a more challenging example because high-conductivity features
appear in many parts of the domain. Numerical results show that our approach
provides an accurate approximation of the solution using a few degrees of freedom and the solution is sparse in an appropriate snapshot space.

Numerical results for the second approach involve solving the Helmholtz
equation in media with two isolated heterogeneous inclusions. We consider a domain with two distinct properties, where plane wave
solutions can provide a good approximation.
In this case, the solution is spanned by only a few plane waves, and is, therefore, sparse
in the space of plane waves with many propagating directions.
However, in general, we do not know which plane wave directions are
dominant and our algorithm identifies these directions by using $l_1$ minimization.

The paper is organized in the following way. In the next section,
we present preliminaries and discuss GMsFEM, coarse and fine grid concepts.
In Section 3, we propose our new construction for the online space.
Section 4 is devoted to numerical results. In Section 5, we present conclusions.

\section{Preliminaries}
\label{prelim_NL}

\begin{figure}[htb]
  \centering
  \includegraphics[width=0.65 \textwidth]{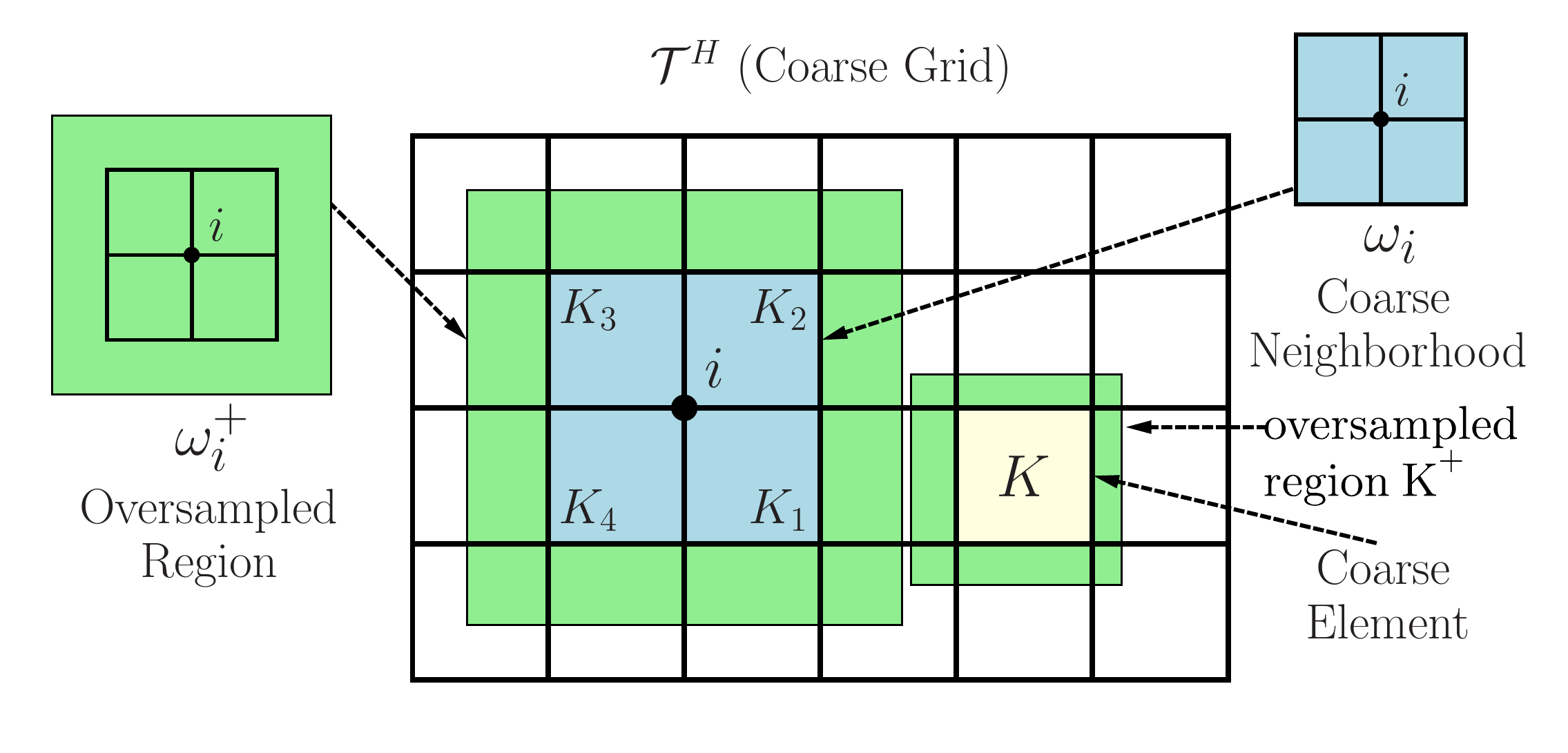}
  \caption{Illustration of a coarse neighborhood and coarse element}
  \label{schematic_ov}
\end{figure}

To discretize (\ref{eq:original_NL}) or (\ref{eq:Helmholtz}),
we let $\mathcal{T}^H$ be a usual conforming partition of the computational domain $D$ into finite elements (triangles, quadrilaterals, tetrahedrals, etc.) and $\mathcal{E}^H$ denotes all the edges in the coarse mesh $\mathcal{T}^H$. We refer to this partition as the coarse grid and assume that each coarse subregion is partitioned into a connected union of fine-grid blocks. The fine-grid partition will be denoted by $\mathcal{T}^h$. We use $\{x_i\}_{i=1}^{N_v}$ (where $N_v$ denotes the number of coarse nodes) to denote the vertices of
the coarse mesh $\mathcal{T}^H$, and define the neighborhood of the node $x_i$ by

\begin{equation} \label{neighborhood}
\omega_i=\bigcup\{ K_j\in\mathcal{T}^H; ~~~ x_i\in \overline{K}_j\}.
\end{equation}
See Figure~\ref{schematic_ov} for an illustration of neighborhoods and elements subordinated to the coarse discretization. We emphasize the use of $\omega_i$ to denote a coarse neighborhood, and $K$ to denote a coarse element throughout the paper.

Next, we briefly outline the global coupling and the role of coarse basis
functions for the respective formulations that we consider. For the discontinuous Galerkin (DG) formulation, we use a coarse element $K$ as the support for basis
functions, and for the continuous Galerkin (CG) formulation, we use $\omega_i$ as the support of basis functions.
In turn, throughout this chapter, we use the notation
\begin{equation} \label{cgordg}
\tau_i = \left\{ \begin{array}{cc}
\omega_i & \text{for} \, \, \,  \text{CG}  \\
K_i & \text{for} \, \, \,  \text{DG}  \\
\end{array}\right.
\end{equation}
when referring to a coarse region where respective local computations are performed (see Figure \ref{schematic_ov}).
To further motivate the coarse basis construction, we offer a brief outline of the global coupling.
In particular, we note that
our approach will employ multiple basis functions per coarse neighborhood.
Both CG and DG solutions will be sought
as $u^{\text{DG/CG}}_{\text{ms}}(x;\mu)=\sum_{i,k} c_{k}^i \psi_{k}^{\tau_i}(x; \mu)$, where $\psi_{k}^{\tau_i}(x; \mu)$ are the basis functions
(without loss of generality, we write basis functions as parameter-dependent).
Once the basis functions are identified, the global coupling is given through the variational form
\begin{equation}
\label{eq:globalG_cg}
a_{\text{DG/CG}}(u^{\text{DG/CG}}_{\text{ms}},v;\mu)=(f,v), \quad \text{for all} \, \, v\in
V_{\text{on}}^{\text{DG/CG}},
\end{equation}
where $V_{\text{on}}^{\text{DG/CG}}$ is used to denote the space formed by those basis functions and $a_{\text{DG/CG}}$ is a bilinear form which will be defined later on. Throughout, for the convenience,
we use the same notations for the discrete
and continuous representations of spatial fields.

\section{Sparse GMsFEM}
\label{cgdgmsfem}

In this section, we will give the detailed constructions of our sparse GMsFEM.
We start with an outline of the approach.

\subsection{Outline}
\label{locbasis}

In this section, we present an outline of the algorithm.
Assume that the snapshot space is
$V_{\text{snap}}^{\tau} = \text{Span}\{\psi_i^{\text{snap}}  \}$
for a generic element $\tau$.
We assume that the solution is sparse in this snapshot space and
consider two approaches (see Figure \ref{overview1} for illustration).
Throughout the paper, we will assume that the solution is sparse
in the snapshot space. We
 will discuss the assumption on the sparsity later in
Section \ref{sec:conclusion}.

\vspace{0.3cm}
\noindent
{\bf General outline of the sparse GMsFEM}:
\begin{itemize}

\item[1.] Coarse grid generation.

\item[2.] Construction of snapshot space, where the solution is sparse.

\item[3.]
\begin{itemize}
\item
{\it First approach.  Local-Sparse Snapshot Subspace Approach}. Seek a subspace of the snapshot space and
construct multiscale basis functions that are sparse in the snapshot space.

\item
{\it Second approach. Sparse Snapshot Space Approach}. Solve for the sparse solution
 in the snapshot space directly within a global formulation.

\end{itemize}

\end{itemize}

In the first approach, we perform local
calculations to identify multiscale basis functions that are
sparse in the snapshot space. Here, we will use approaches similar
to sparse POD. Then, the global problem is solved in the space
of multiscale basis functions. The resulting solution is sparse.

In the second approach,
we apply directly sparse solution techniques and find the solution that is sparse
in the snapshot space. This approach is more expensive
since it uses a large snapshot space. However, in some examples,
we can not identify local basis functions in the offline stage and
such approaches
give sparse solutions in the online stage.

\subsection{First approach.  Local-Sparse Snapshot Subspace Approach}

We first give a general idea of this approach. We consider
a local snapshot space $V_{\text{snap}}^{\tau}=\text{Span}\{\psi_i^{\text{snap}}  \}$.
In the local snapshot space, we seek multiscale basis functions $\{\psi_i^{\text{on}}\}$
that are sparse in the local snapshot space and which have
smallest energies (similar to sparse POD). For example,
following to \cite{schaeffer2013sparse}, we can consider
\[
\min\limits_{\Psi\in \mathbb{R}^{n\times M_{\text{on}}^{\tau}}}{1\over \nu}\norm{\Psi}_{1}+\text{Tr}\langle \Psi^{T}A_{s}(\mu)\Psi\rangle,\text{ s.t. }\Psi^T\Psi=\text{I},
\]
where $n$ is the dimension of the fine-scale space, $M_{\text{on}}^{\tau}$ is the number of online basis functions
and $A_s(\mu)$ is the stiffness matrix formed in the snapshot space.
Our proposed local approach avoids expensive direct local
eigenvalue calculations and uses randomized snapshots.
This approach will be
 used for problems where one knows a local basis construction principle.
The latter typically uses
some local generalized eigenvalue problems in the
snapshot spaces \cite{schaeffer2013sparse}.
Once multiscale basis functions (that are sparse in the snapshot space)
are constructed, we solve the problem on a coarse grid.

Here, we will consider the parameter-dependent problem
(\ref{eq:original_NL}).
We can consider the computation of the parameter-dependent coarse space
as an online procedure.
In the latter, our objective will be solving
many problems for a given value of
the parameter with different boundary conditions and the right hand sides.

\subsubsection{Snapshot space}

We first construct a snapshot space $V_{\text{snap}}^{\tau}$
(for a generic $\tau$).
 Construction of the snapshot space involves solving the local problems for various choices of input parameters, and we describe it below.
We generate  snapshots using random boundary conditions by
solving a small number of local problems imposed with random boundary conditions,
\begin{equation}
\label{eq:random bc}
 \begin{aligned}
 -\mbox{div}(\kappa (x; \mu_j)\nabla  \psi_{l,j}^{\tau, \text{rsnap}}) &=0\ \ \text{in } \tau^+\\
  \psi_{l,j}^{\tau, \text{rsnap}}&=r_{l} \text{ on } \partial\tau^+,
 \end{aligned}
 \end{equation}
 where $r_{l}$ are independent identically distributed (i.i.d.) standard Gaussian random vectors on
 the fine-grid nodes of the boundary, $l=1, \cdots,L$. $\mu_j$ ($j=1,\dots,J$) is a specified set of fixed parameter values, and $J$ denotes the number of parameter values we choose. Here, $\tau^+$ is an oversampled region shown in the Figure \ref{schematic_ov} as $\omega_i^{+}$ or $K^+$ for conforming Galerkin formulation or discontinuous Galerkin formulation.

 The space generated by $\psi_{l,
 j}^{ \tau, \text{rsnap}}$ is a subspace of the space generated by all local snapshots
 $\Psi_{k, j}^{\tau, \text{snap}}$, $k=1,\cdots, N$, and $N$ denotes the number of boundary nodes.
 Denote $\Psi_{j}^{\tau, \text{rsnap}}=[\psi_{1,j}^{\tau, \text{rsnap}}, \cdots, \psi_{L,j}^{\tau, \text{rsnap}}]$ and $\Psi_{j}^{\tau, \text{snap}}=[\psi_{1,j}^{\tau, \text{snap}}, \cdots, \psi_{N,j}^{\tau, \text{snap}}]$.
 Therefore, for each parameter $\mu_j$, there exists a randomized matrix $\mathcal{R}$
 with rows composed by the random boundary vectors $r_{l}$ (as shown in Equation \eqref{eq:random bc}), such that,
\begin{align}\label{eqn:random_snapshots}
 \Psi_{j}^{\tau, \text{rsnap}}=\mathcal{R}\Psi_{j}^{\tau, \text{snap}}.
 \end{align}
Now, we are ready to present the local snapshot space as follows,
$$
V_{\text{snap}}^{\tau} = \text{Span}\{ \psi_{l,j}^{ \tau, \text{rsnap}}:~~~1\leq j \leq J ~~ \text{and} ~~ 1\leq l \leq L \},
$$
for each coarse subdomain $\tau$.
\begin{remark}
Note that we impose the same random vectors for the local snapshot calculation in Equation \eqref{eq:random bc} for $\mu_j$, $j=1,\cdots, J$ in order to obtain a sparse online space.
\end{remark}

\subsubsection{Sparse local space calculations}

For a given input parameter $\mu$, we next construct the associated
coarse space
$V^{\tau}_{\text{on}}(\mu)$ \emph{for each} $\mu$ value on each coarse subdomain $\tau$.
In principle, we want this to be a small dimensional subspace of the snapshot space for computational efficiency.
The  coarse space will be used within the finite element
framework to solve the original global problem, where a continuous or discontinuous Galerkin coupling of the multiscale basis functions is used to compute the global solution. In particular, we seek a subspace of the snapshot space $V_{\text{snap}}^{\tau}$ such that it can approximate any element of the snapshot space in
an appropriate sense. For the convenience of the presentation, we
denote $\Psi_{\text{snap}}^{\tau}=[\psi_{1}^{\tau, \text{snap}},\cdots,\psi_{J}^{\tau, \text{snap}}]$.
Similar as the generation of local snapshot space, we obtain the local sparse online basis via a local problem solved by $l_1$ optimization in the space of the corresponding local snapshot space. Here, we will use a smaller subspace of $V_{\text{snap}}^{\tau}$ as the test space constructed through multiplication of a random matrix $T_{\text{random}}^{\tau}$. Denote $T_{\text{random}}^{\tau}$ as a matrix of size $(L\times J) \text{by }q$ ($q<< (L\times J)$) with rows of i.i.d standard Gaussian random vectors. Then the test space is defined as $\Psi_{\text{snap}}^{\tau}T_{\text{random}}^{\tau}$.

Specifically, the local problem is arranged as follows.
Find $U_{l}$, such that, $\psi_{l}^{\tau, \text{on}}=\Psi_{\text{snap}}^{\tau}U_{l}$,  and
\begin{align}
\label{eq:l1opt}
U_{l}=\text{argmin} {1\over \nu}\norm{U}_{1} \text{ subject to } A_c(\mu)U_{l}=F_{l}.
\end{align}
Here,
 $A_c(\mu)=({\Psi_{\text{snap}}^{\tau}T_{\text{random}}^{\tau}})^T A(\mu) \Psi_{\text{snap}}^{\tau}$ and $F_{l} = ({\Psi_{\text{snap}}^{\tau}T_{\text{random}}^{\tau}})^T R_{l}$ with $A(\mu)$ being the local stiffness matrix and $R_{l}$ being the right hand side for the local problem with Dirichlet boundary condition $r_{l}$. Namely, we are solving the following local problems in $V_{\text{snap}}^{\tau}$ with the $l_1$ minimized coefficient vector $U_l$ for the testing space equal to  $\Psi_{\text{snap}}^{\tau}T_{\text{random}}^{\tau}$,
\begin{equation}
\label{eq:online}
\left\{
 \begin{aligned}
 -\mbox{div}(\kappa (x; \mu)\nabla  \psi_{l}^{\tau^+, \text{on}})=0\ \ \text{in}\  \tau^+\\
  \psi_{l}^{\tau^+, \text{on}}=r_{l} \text{ on } \partial\tau^+.
 \end{aligned}
 \right.
 \end{equation}
 Later, we will briefly introduce the algorithm to solve Equation \eqref{eq:l1opt}.
Note that we impose the same random vectors as the boundary conditions in Equation \eqref{eq:random bc} to guarantee the sparse solution in Equation \eqref{eq:online}.

We can obtain the local online snapshot functions on the target
domain $\tau$ by restricting the solution of the local problem, $\psi_{l}^{\tau^+, \text{on}}$
to $\tau$ (which is denoted by $\psi_{l}^{\tau, \text{on}}$).
Now we are ready to present the local online snapshot space as follows,
$$
V_{\text{on}}^{\tau} = \text{Span}\{ \psi_{l}^{\tau,  \text{on}}:~~~1\leq l \leq L \},
$$
for each coarse subdomain $\tau$. Then we denote $\Psi_{\text{on}}^{\tau}=[\psi_{1}^{\tau, \text{on}}, \cdots, \psi_{L}^{\tau, \text{on}}]$.

Next we will select the dominant modes from $V_{\text{on}}^{\tau}$ by the following eigenvalue problem,
\begin{equation} \label{oneig}
A^{\tau, \text{on}}(\mu) z_k^{\tau, \text{on}} = \lambda_k^{\tau, \text{on}} S^{\tau, \text{on}}(\mu) z_k^{\tau, \text{on}},
\end{equation}
where
\begin{equation*}
 \displaystyle A^{\tau, \text{on}}(\mu) = [a^{\text{on}}(\mu)_{mn}] = [\int_\tau \kappa(x; \mu)
 \nabla \psi_m^{\tau, \text{on}} \cdot \nabla \psi_n^{\tau, \text{on}}] = {\Psi_{\text{on}}^{\tau}}^T A(\mu) \Psi_{\text{on}}^{\tau}
 \end{equation*}
\begin{equation*}
 \displaystyle S^{\tau, \text{on}}(\mu) = [s^{\text{on}}(\mu)_{mn}] = [\int_{\tau} {\kappa}(x; \mu)
 \psi_m^{\tau, \text{on}} \psi_n^{\tau, \text{on}}] = {\Psi_{\text{on}}^{\tau}}^T S(\mu) \Psi_{\text{on}}^{\tau},
 \end{equation*}
and $\kappa(x; \mu)$ is now parameter dependent.
To generate the coarse space we then choose the smallest $M_{\text{on}}^{\tau}$ eigenvalues from Equation \eqref{oneig} and form the corresponding eigenvectors in $R_{\text{on}}^{\tau}$ by setting
$\phi_k^{\tau,\text{on}} = \sum\limits_{j=1}^{L}\psi_{j}^{\tau,\text{on}} z_{k,j}^{\tau,\text{on}}$
(for $k=1,\ldots, M_{\text{on}}^{\tau}$), where $z_{k,j}^{\tau,\text{on}}$ are the coordinates of
the vector $\phi_{k}^{\tau,\text{on}}$.

%


Above, we presented an algorithm for solving an eigenvalue problem
in the space which has a sparse representation in the snapshot space.
One can attempt to obtain sparse spectral basis from Equation \eqref{eq:online} using $l_1$-minimization method in the local snapshot space $V_{\text{snap}}^{\tau}$ (\cite{hale2008fixed, yin2008bregman}).
Here, we can use the algorithm proposed in \cite{ozolins2014compressed}.
The general basis pursuit problem is as follows
\begin{align}
\label{eq:l1prb}
\min\limits_{x\in \mathbb{R}^n}\norm{x}_{1}+\nu\norm{Cx-f(x)},
\end{align}
where $f\in\mathbb{R}^m$, $C\in\mathbb{R}^{m\times n}$, and $m<\!\!<n$.
We refer to \cite{yin2008bregman} for the Bregman algorithm to solve \eqref{eq:l1prb}.
Instead of solving Equation \eqref{eq:online} in the online stage, we can generate the sparse spectral basis directly following the algorithm in
\cite{ozolins2014compressed}. The first $M_{\text{on}}^{\tau}$ sparse spectral basis can be solved by
\begin{align}
\label{eq:sparse_eig}
\min\limits_{\Psi\in \mathbb{R}^{n\times M_{\text{on}}^{\tau}}}{1\over \nu}\norm{\Psi}_{1}+\text{Tr}\langle \Psi^{T}A_{c}(\mu)\Psi\rangle,\text{ s.t. }\Psi^T\Psi=\text{I}.
\end{align}

\subsubsection{Global coupling}
\label{globcoupling}

The multiscale basis functions constructed above can be coupled via DG or CG
formulation. Below, we present DG approach (similar results are
observed when CG approach is used).
One  uses the discontinuous Galerkin (DG) approach
(see also \cite{riviere2008discontinuous,ABCM_unified_2002})
to couple multiscale basis functions. This may avoid the use of
the partition of unity functions;
however, a global formulation needs to be chosen carefully.
The global formulation is given by
\begin{equation}
a_{\text{DG}}(u_{H}^{\text{DG}},v)=(f,v),\quad\forall v\in V_{\text{on}},\label{eq:ipdg}
\end{equation}
where the bilinear form $a^{\text{DG}}$ is defined as
\begin{equation}
a_{\text{DG}}(u,v)=a_{H}(u,v)-\sum_{E\in\mathcal{E}^{H}}\int_{E}\Big(\average{{\kappa}\nabla{u}\cdot{n}_{E}}\jump{v}+\average{{\kappa}\nabla{v}\cdot{n}_{E}}\jump{u}\Big)+\sum_{E\in\mathcal{E}^{H}}\frac{\gamma}{h}\int_{E}\overline{\kappa}\jump{u} \jump{v} \label{eq:bilinear-ipdg}
\end{equation}
with
\begin{equation}
a_{H}({u},{v})=\sum_{K\in\mathcal{T}_{H}}a_{H}^{K}(u,v),\quad a_{H}^{K}(u,v)=\int_{K}\kappa\nabla u\cdot\nabla v,
\end{equation}
where $\gamma>0$ is a penalty parameter, ${n}_{E}$ is a fixed unit
normal vector defined on the coarse edge $E \in \mathcal{E}^H$.
Note that, in (\ref{eq:bilinear-ipdg}),
the average and the jump operators are defined in the classical way.
Specifically, consider an interior coarse edge $E\in\mathcal{E}^{H}$
and let $K^{+}$ and $K^{-}$ be the two coarse grid blocks sharing
the edge $E$. For a piecewise smooth function $G$, we define
\[
\average{G}=\frac{1}{2}(G^{+}+G^{-}),\quad\quad\jump{G}=G^{+}-G^{-},\quad\quad\text{ on }\, E,
\]
where $G^{+}=G|_{K^{+}}$ and $G^{-}=G|_{K^{-}}$ and we assume that
the normal vector ${n}_{E}$ is pointing from $K^{+}$ to $K^{-}$.
Moreover, on the edge $E$, we define $\overline{\kappa} = (\kappa_{K^+}+\kappa_{K^-})/2$
where $\kappa_{K^{\pm}}$ is the maximum value of $\kappa$ over $K^{\pm}$.
For a coarse edge $E$ lying on the boundary $\partial D$, we define
\[
\average{G}=\jump{G}=G,\quad \text{ and }\quad \overline{\kappa} = \kappa_{K} \quad\quad\text{ on }\, E,
\]
where we always assume that ${n}_{E}$ is pointing outside of $D$.
We note that the DG coupling (\ref{eq:ipdg})
is the classical interior penalty discontinuous Galerkin (IPDG) method
\cite{riviere2008discontinuous}
with our multiscale basis functions as the approximation space.

We can obtain the discontinuous Galerkin spectral multiscale space as
\begin{equation} \label{dgspace}
V_{\text{on}}^{\text{DG}}(\mu) = \text{Span} \{ \phi_{k}^{K,\text{on}}: \,  \, \, 1 \leq k \leq M_{\text{on}}^{K}, K \in \mathcal{T}^H  \}.
\end{equation}

We can obtain an operator matrix constructed by the basis functions of $V_{\text{on}}^{\text{DG}}(\mu)$. We denote the matrix as $\Phi_{0}$ where $\Phi_0 = \left[ \phi_1^{\text{DG}} , \ldots, \phi_{N_c}^{\text{DG}} \right]$. Recall that $N_c$ denotes the total number of coarse basis functions.
Solving the problem \eqref{eq:original_NL} in the coarse space $V_{\text{on}}^{\text{DG}}(\mu)$ using the DG formulation described in Equation \eqref{eq:ipdg} is equivalent to seeking
$u^{\text{DG}}_{\text{ms}}(x; \mu) = \sum_i c_i \phi_i^{\text{DG}}(x; \mu) \in V_{\text{on}}^{\text{DG}}$ such that
\begin{equation} \label{dgvarform}
a^{\text{DG}}(u_{\text{ms}}^{\text{DG}}, v; \mu) = (f, v) \quad \text{for all} \,\,\, v \in V_{\text{on}}^{\text{DG}},
\end{equation}
where
$ \displaystyle a^{\text{DG}}(u, v; \mu) $ and $f( v)$ are defined in Equation \eqref{eq:bilinear-ipdg}.
We can obtain a coarse  system
\begin{equation}
A_0 U_0^{\text{DG}} = F_0,
\end{equation}
where $U^{\text{DG}}_0$ denotes the discrete coarse DG solution, and
\begin{equation*}
A_0(\mu) = R_0^T A(\mu) R_0 \quad \text{and} \quad F_0 = R_0^T F,
\end{equation*}
where $A(\mu)$ and $F$ are the standard, fine-scale stiffness matrix and forcing vector corresponding to the form in Equation \eqref{dgvarform}. After solving the coarse system, we can use the operator matrix $R_0$ to obtain the fine-scale solution in the form of $R_0U_0^{\text{DG}}$.


\subsubsection{Computational cost}

In this section, we discuss the computational cost. For this, we assume
that we have chosen $J$ parameters,
$\mu_1,..., \mu_J$,  and
$L$ boundary conditions $r_1,..., r_L$,
for constructing the snapshot space.
Then, the cost for snapshot calculations will be the same as solving
$J \times L$ local problems for randomized snapshots.
Next, we compute the cost of solving
$L_{\text{on}}$  online randomized snapshots. Each online snapshot
calculation requires solving $l_1$ minimization with a constraint
involving $q \times (L \times J)$ matrix (see (\ref{eq:l1opt})).
The cost of eigenvalue computation with
$L_{\text{on}}$ snapshots is considered to be small as it involves a
small eigenvalue problem of the size $L_{\text{on}}\times L_{\text{on}}$.
The online cost is mainly due to solving (1) solving
$L_{\text{on}}$  online randomized snapshots (2) solving a global problem
on a coarse grid. The cost of solving a global problem is small
if the solution has a sparse representation. As for the cost of solving
online randomized snapshots, this will be small compared to solving local
problems in the online stage if local problems have a high resolution.
 Moreover, as we pointed out
earlier
that the proposed approach allows a fast assembly of the stiffness
matrix in the online space since it involves a few elements of the
snapshot space. It also avoids using all snapshot vectors,
 which can result to a large dimensional representation
of the online multiscale basis functions.


\subsection{Second approach. Sparse Snapshot Subspace Approach}

In this approach, we will use
an appropriate snapshot space to compute the sparse solution directly.
Again, we consider
a local snapshot space $V_{\text{snap}}^{\tau}=\text{Span}\{\psi_i^{\tau, \text{snap}}  \}$.
In some applications, we may not be able to reduce the dimension of the multiscale space locally. In this case, we can use sparsity techniques directly
in the global snapshot space to compute the solution.
The procedure can be more expensive; however, can yield more accurate solutions.
More precisely, we seek the solution in the global snapshot space of
$V_{\text{snap}} = \text{Span}\{ \psi_i^{\text{snap}} \}
$
using $l_1$ minimization with testing space, $V_\text{test}$, spanned by the random combination of snapshot basis functions. This is equivalent to find $u_\text{ms}=\sum_i U^\text{ms}_i \psi^\text{snap}_i\in V_\text{snap} $ where
\begin{align}
\label{eq:l1_global}
U^\text{ms}=\text{argmin} \norm{U}_{1} \text{ subject to } a_\text{DG}(\sum_i U_i \psi^\text{snap}_i,v)=(f,v), \;\forall v\in V_\text{test}.
\end{align}

In the following section, we will take the Helmholtz problem as an example and discuss the procedure of using this approach to compute a sparse solution. More precisely, we consider the following problem: find $u$ such that
\[
 -\nabla \cdot( \kappa(x) \nabla u) - \Omega^2 n(x)u = f, \quad \text{in} \quad D
\]
with the Dirichlet boundary condition $u|_{\partial D }=g$.

\subsubsection{Snapshot and test space}

In this section, we present the construction of the snapshot space $V_\text{snap}$ and the test space $V_\text{test}$. Since we are solving the Helmholtz equation with a fixed frequency $\Omega$, we can assume the solution, $u$, can be written as a linear combination of plane waves, namely, $u = \sum_k \beta_k e^{i \Omega k\cdot x}$. Therefore, we consider our snapshot basis to be some plane waves in each coarse block $K\in \mathcal{T}^H$, that is,
\[
V_\text{snap} = \text{Span}\{ \psi_{m,j}:\,1\leq m \leq N_d \text{ and } 1\leq j \leq M\},
\]
where
\begin{equation}
\label{eq:SparseTrig}
\Psi_{m,j}(x) = \left\{ \begin{array}{cc}
 e^{i \Omega k_m\cdot x} & \text{for} \, \, \,  x\in K_j  \\
0 & \text{otherwise}  \\
\end{array}\right.
\end{equation}
with $k_m= (\sin(\pi m/N_d),\cos(\pi m/N_d)) $
, $N_d$ as the number of propagating directions and $M$ the number of coarse blocks.
We note that the plane wave basis is not new in solving the Helmholtz equation
and it was used in a number of papers (see \cite{huttunen2007use, colton1985novel,colton1987numerical,tezaur2006three, hiptmair2011plane} and
the references therein).

Next, we will show the construction of the test space. We consider $\{ r^{(l)} \}^{N_t}_{l=1}$ as a collection of i.i.d. standard Gaussian random vectors with $N_t<<N_d$. Then the testing space $V_\text{test}$  is defined by
\[
V_\text{test}=\text{Span} \{ \phi_{l,j}: \phi_{l,j} =\sum_{m=1}^{N_d} r^{(l)}_{m}\psi_{m,j}, \,1\leq l \leq N_t \text{ and } 1\leq j \leq M\ \}.
\]
Since $N_t<<N_d$, we have a test space with dimension much smaller the snapshot space ($\text{dim}(V_\text{test}) = N_t M<<N_d M=\text{dim}(V_\text{snap} )$).

\subsubsection{Sparse solution in the snapshot space}

After constructing the snapshot and test space, we can couple the system globally by IPDG method. That is, we find $u_\text{ms}\in V_{\text{snap}}$ such that
\begin{equation}
\label{eq:dg_app2}
a_\text{DG}(u_\text{ms},v) = (f,v),\,\quad \forall v\in V_\text{test},
\end{equation}
where \begin{equation}
a_{\text{DG}}(u,v)=a_{H}(u,v)-\sum_{E\in\mathcal{E}^{H}}\int_{E}\Big(\average{{\kappa}\nabla{u}\cdot{n}_{E}}\jump{v}+\average{{\kappa}\nabla{v}\cdot{n}_{E}}\jump{u}\Big)+\sum_{E\in\mathcal{E}^{H}}\frac{\gamma}{h}\int_{E}\overline{\kappa}\jump{u} \jump{v} \label{eq:bilinear-ipdg2}
\end{equation}
with
\begin{equation}
a_{H}({u},{v})=\sum_{K\in\mathcal{T}_{H}}a_{H}^{K}(u,v),\quad a_{H}^{K}(u,v)=\int_{K}\kappa\nabla u\cdot\nabla v - \Omega^2 \int_{K}n(x) uv,
\end{equation}
where $\gamma>0$ is a penalty parameter, ${n}_{E}$ is a fixed unit
normal vector defined on the coarse edge $E \in \mathcal{E}^H$. Moreover, $\{ \cdot \}$ and $[\cdot]$ are the average and jump operators defined before.

As the dimension of test space is smaller that the dimension of snapshot space, the Equation \eqref{eq:dg_app2} does not have a unique solution. To seek for a sparse solution, we will solve a $l_1$ minimization problem subject to Equation \eqref{eq:dg_app2}, more precisely, we will find $u_\text{ms}=\sum_{l,j} U^\text{ms}_{l,j}\phi_{l,j} \in V_\text{snap}$ such that
\begin{equation}
\label{eq:l1-helmholtz}
U^\text{ms} = \text{argmin} \{ \|U \|_{\l_1}\} \text{ subject to } a_{DG}(U^\text{ms}_{l,j}\phi_{l,j},v)=(f,v),\quad\;\forall v\in V_\text{test}.
\end{equation}

\subsubsection{Cost of computations}

In this section, we discuss the computational cost associated
with the second approach.
In this example,
the cost for snapshot (plane waves)
calculations is cheap as they are analytically
described.
In addition, the linear system in (\ref{eq:l1-helmholtz}) has dimension $N_t M \times N_d M$,
which is a highly under-determined system as $N_t << N_d$.
Thus, one can solve the $l_1$ minimization problem
(\ref{eq:l1-helmholtz}) efficiently
by using, for example, the Bregman's method \cite{yin2008bregman}.
Moreover, if an adaptivity can be used and the problem requires
very few snapshot vectors or very few test vectors
in many regions excepts a few coarse regions,
this will increase the efficiency of the proposed approach.
In this case, if we denote by $N_t^{(i)}$ and $N_d^{(i)}$
the number of test and snapshot vectors
in the region $i$ (using adaptivity), then the
 linear system in (\ref{eq:l1-helmholtz}) has dimension $\left(\sum_{i=1}^M
N_t^{(i)} \right) \times \left(\sum_{i=1}^M N_d^{(i)} \right)$.
Thus, if $\sum_{i=1}^M N_t^{(i)}$  or $\sum_{i=1}^M N_d^{(i)}$ is not large,
one can gain computational efficiency.
We note that this approach is more expensive compared to the first approach,
where we perform the sparsity calculations at the local coarse-grid level.

\section{Numerical results}

\subsection{First Approach. Local-Sparse Snapshot Subspace Approach}

In this section, we will present some numerical examples by using the first approach to compare the sparse multiscale solution. We consider the domain $D = [0,1]^2$. The coarse mesh size $H$ is $1/10$ and each coarse grid block is subdivided into a $10\times 10$ grid, therefore, the fine mesh size $h=1/100$;

\subsubsection*{Example 1}

\begin{figure}\centering
 \subfigure[$\kappa_1(x)$]{\label{fig:permi}
    \includegraphics[width = 0.45\textwidth, keepaspectratio = true]{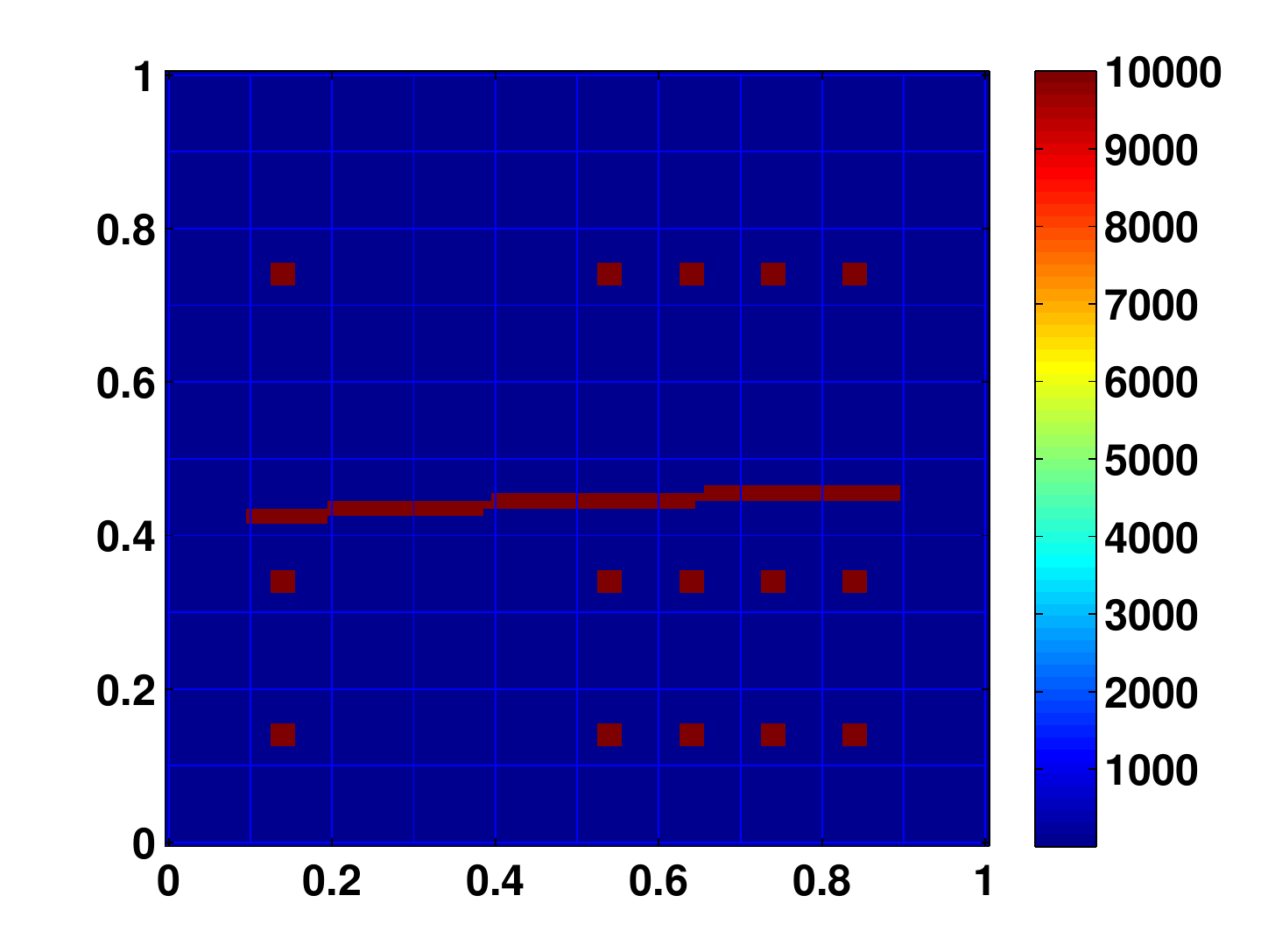}
   }
  \subfigure[$\kappa_2(x)$]{\label{fig:permii}
     \includegraphics[width = 0.45\textwidth, keepaspectratio = true]{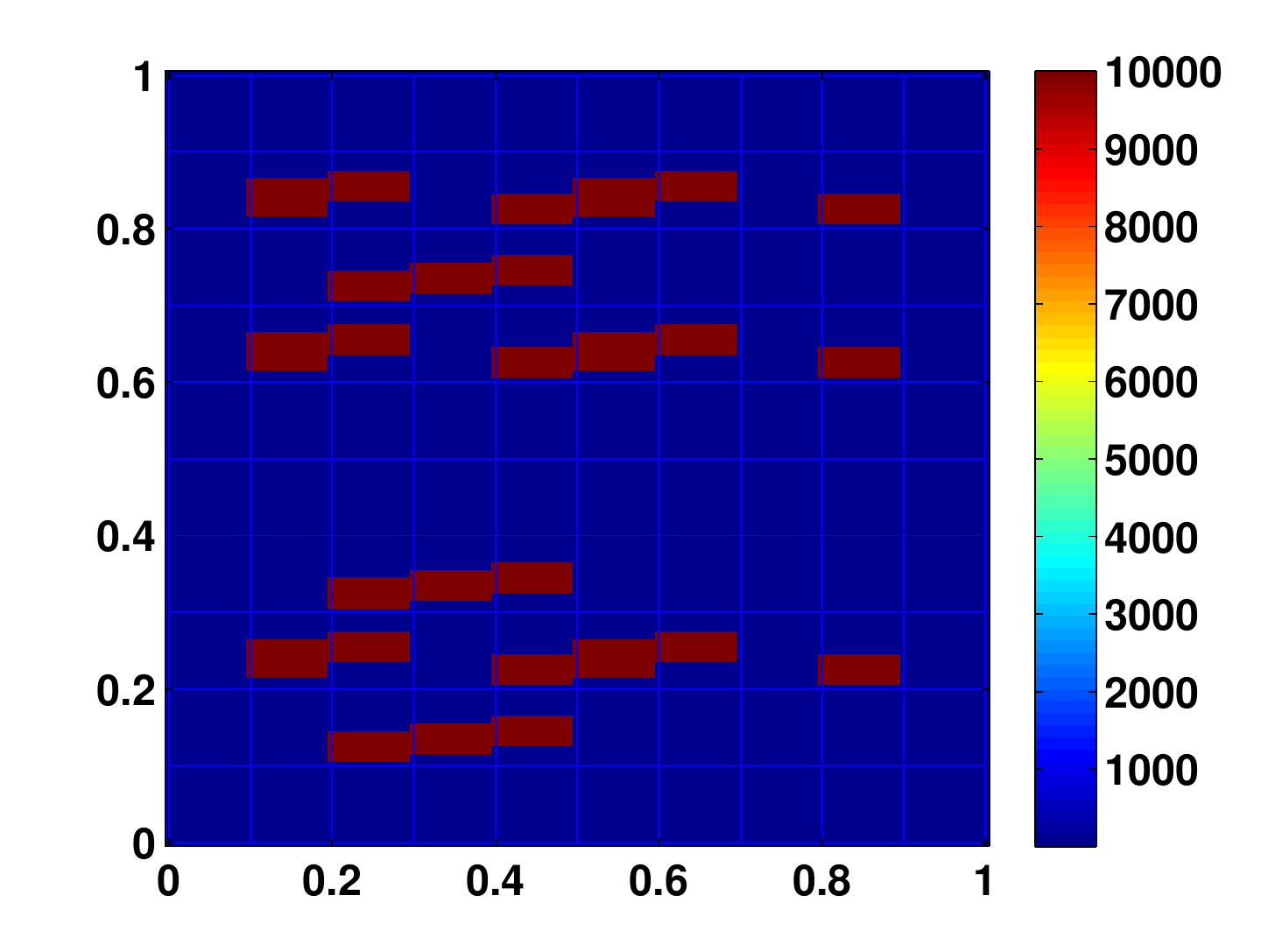}
  }
 \caption{Decomposition of permeability field}\label{fig:decofperm}
\end{figure}

{\bf Setup.}
In our first example, we will consider the source function $f=1$ and the medium parameter $\kappa(\mu) = (1-\mu)\kappa_1 + \mu\kappa_2$,
 where $\kappa_1$ and $\kappa_2$ are shown in Figure \ref{fig:decofperm}. We choose the offline values of
$\mu$, $\mu_i=0.2,0.4,0.6,0.8$,
 for computing the online snapshot space as discussed above.

\begin{figure}[ht]\centering
\includegraphics[scale=0.4]{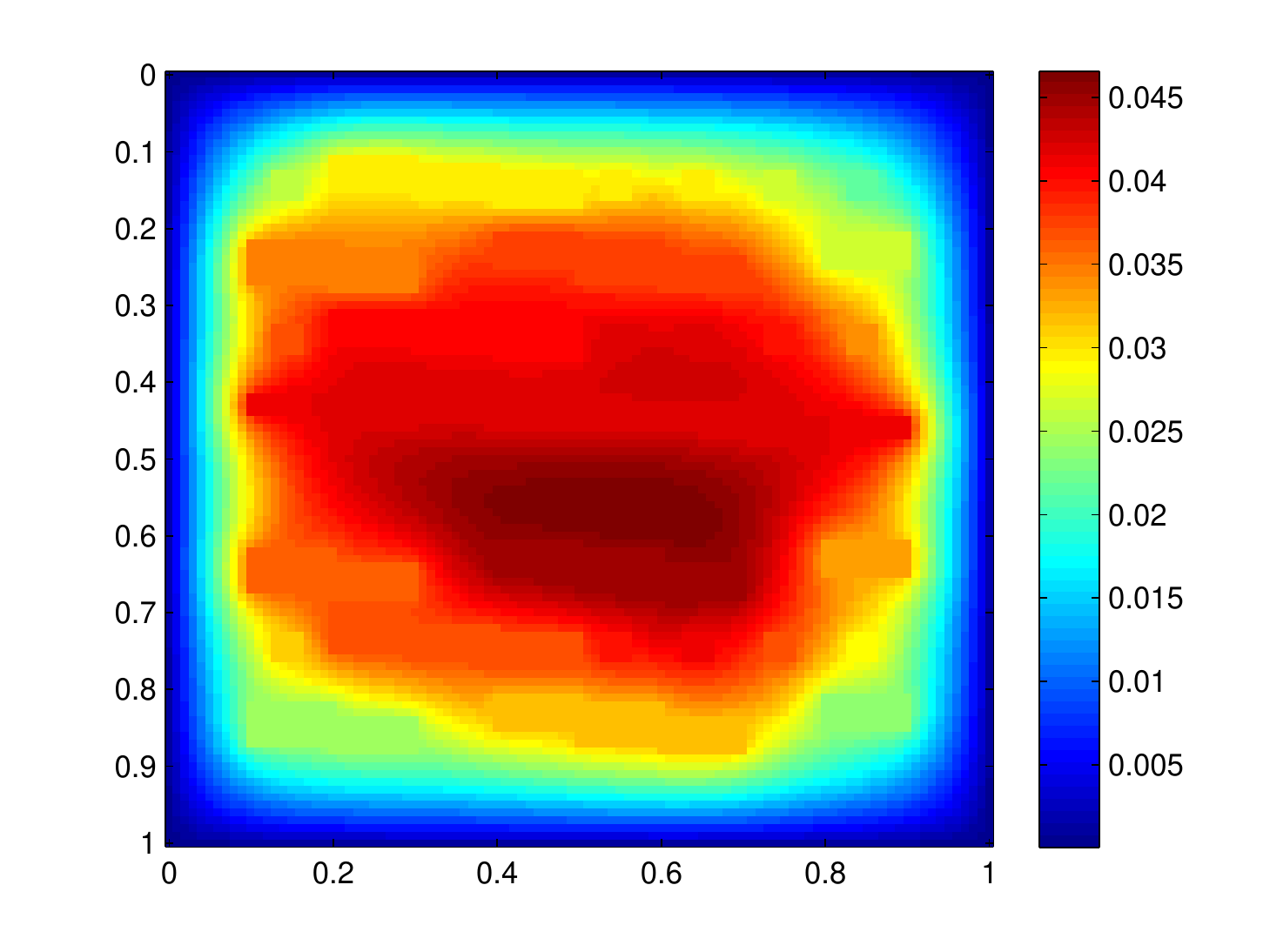}\includegraphics[scale=0.4]{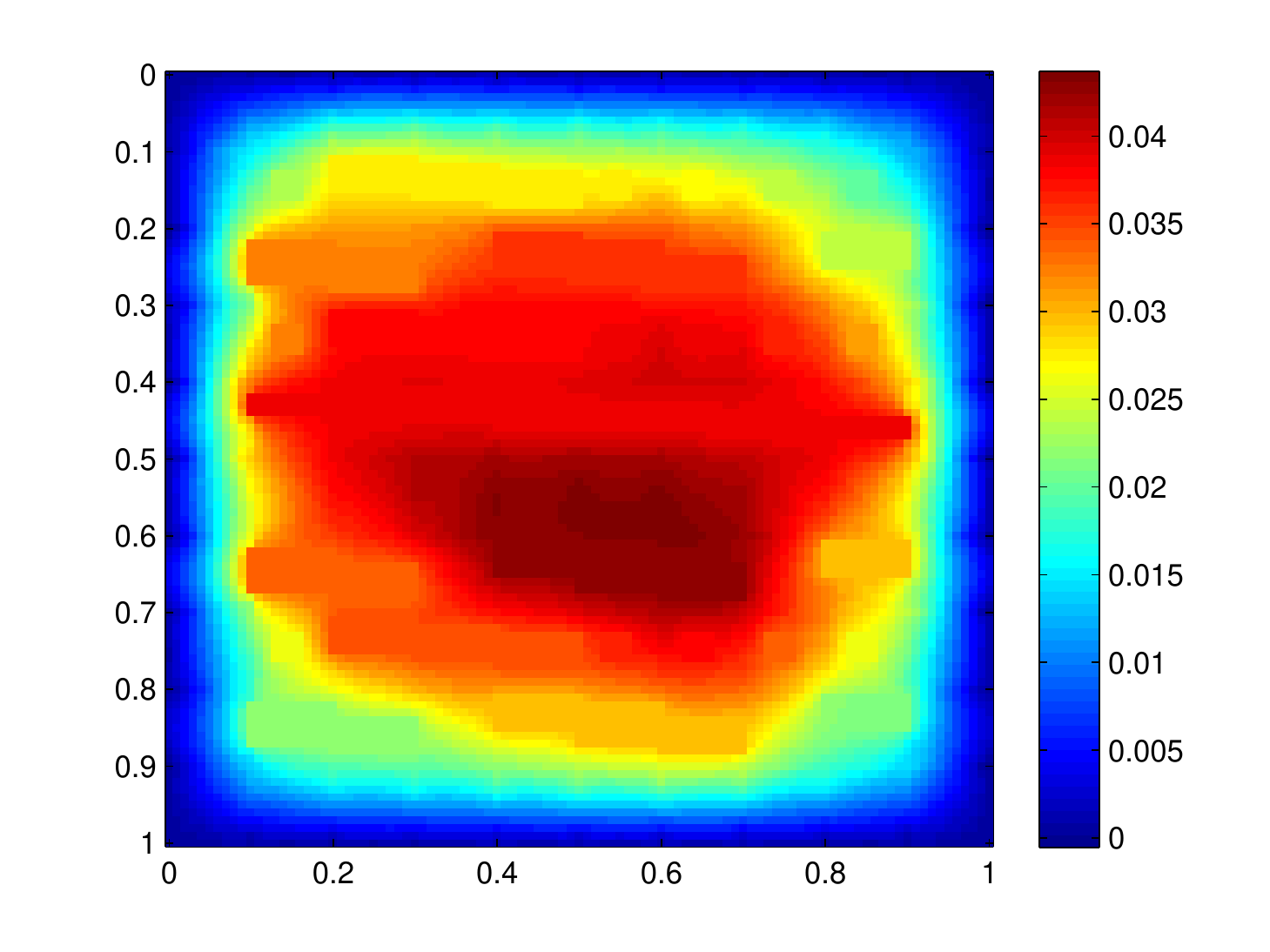}

\protect\caption{Left: Fine grid solution Right: Numerical solution ($8$ basis).}
\label{fig:sol_case1}
\end{figure}

{\bf Discussions of numerical results.}
In Table \ref{table:sparse_DG Harmonic}, we show the convergence history of our method for $\mu=0.5$,
where we define$\| u\|_{H^1_{\kappa}(D)}^2 = \int_D \kappa |\nabla u|^2$.
The fine-grid solution and the numerical solution are shown in Figure \ref{fig:sol_case1}.
First, we note
that there is an irreducible error due to the use of the snapshot
space, which consists of harmonic functions. This error is of order of
the coarse mesh size and this is the reason, the error decay is slow (below
$10$ \%) as
we increase the dimension. In these problems,
because our selected $\mu=0.5$ is near to $\mu=0.4$ and $\mu=0.6$,
we observe that the sparsity is close to $50$ \%, i.e., we only use
snapshots corresponding to nearby values of $\mu$.
We observe that when the snapshot space dimension
is $9600$ (i.e., $24$ randomized solutions per coarse block and per
each value of $\mu_i$), the nonzero coefficients in the expansion of
basis functions (over the whole domain) in terms of $9600$ snapshot
vectors are $4850$. The optimal expansion for this case
corresponds when $\mu=0.5$
is selected for snapshot construction and in this case, the number of
nonzero coefficients is $2400$ ($24$ per coarse region).
We note that if we consider small dimensional online spaces and
the full snapshot space, then the sparsity is very small.
For example, if we use $12$ randomized solutions per coarse block
and per each value of $\mu_i$ for identifying multiscale basis functions
per each coarse region, then, we will be using only $1/2$ of the snapshot
vectors and thus, the sparsity (the number of nonzero coefficients
of the solution in the snapshot space) will be $25$ \%.
As we observe
that our numerical examples identify appropriate sparsity of the solution
space. We expect a more significant gain in the sparsity if more parameter
values are used.


{\bf Why to expect a sparsity.} Next, we briefly describe why to expect
a sparsity in this problem. Because
the snapshot space consists of
local problems corresponding to
multiple values of $\mu$, we expect that for an online value of $\mu$,
we will have a local (in $K$) coefficient $\kappa(x;\mu)$, which
is similar to one of the snapshot solutions. Thus, it is more advantageous
to use $l_1$ minimization techniques, which will select a small dimensional
subspace of the snapshot space corresponding to the coefficient that is close
to $\kappa(x;\mu)$ with the online value of $\mu$. Such situations may occur
in various applications. Moreover, in these examples, it is more advantageous
to use local spectral decomposition and avoid a large-scale $l_1$ minimization
problem.


\begin{table}[htb!]
\centering
\caption{Convergence history of the DGMsFEM using oversampling Harmonic basis. The fine-scale dimension is 12100.
The full snapshot space dimension is 9600.}
 \label{table:sparse_DG Harmonic}
\begin{tabular}{|c|c|c|c|c|}
\hline
\multirow{2}{*}{$\text{dim}(V_{\text{on}})$}  &
\multicolumn{2}{c|}{  $\|u-u_{\text{ms}} \|$ (\%) }  \\
\cline{2-3} {}&
$\hspace*{0.8cm}   L^{2}(D)   \hspace*{0.8cm}$ &
$\hspace*{0.8cm}   H^{1}_\kappa(D)  \hspace*{0.8cm}$
\\
\hline\hline
       $400$     &   $15.05$    & $31.84$  \\
\hline
      $600$    &    $2.89$    & $13.71$ \\
\hline
      $800$    &    $1.22$    & $10.20$ \\
\hline
     $1000$   &   $1.12$       & $9.83$  \\
\hline
     $\text{dim}(V_{\text{snap}})$   &   $1.07$       & $9.59$  \\
\hline
\end{tabular}
\end{table}


\subsubsection*{Example 2}

\begin{figure}[ht]\centering
\includegraphics[scale=0.4]{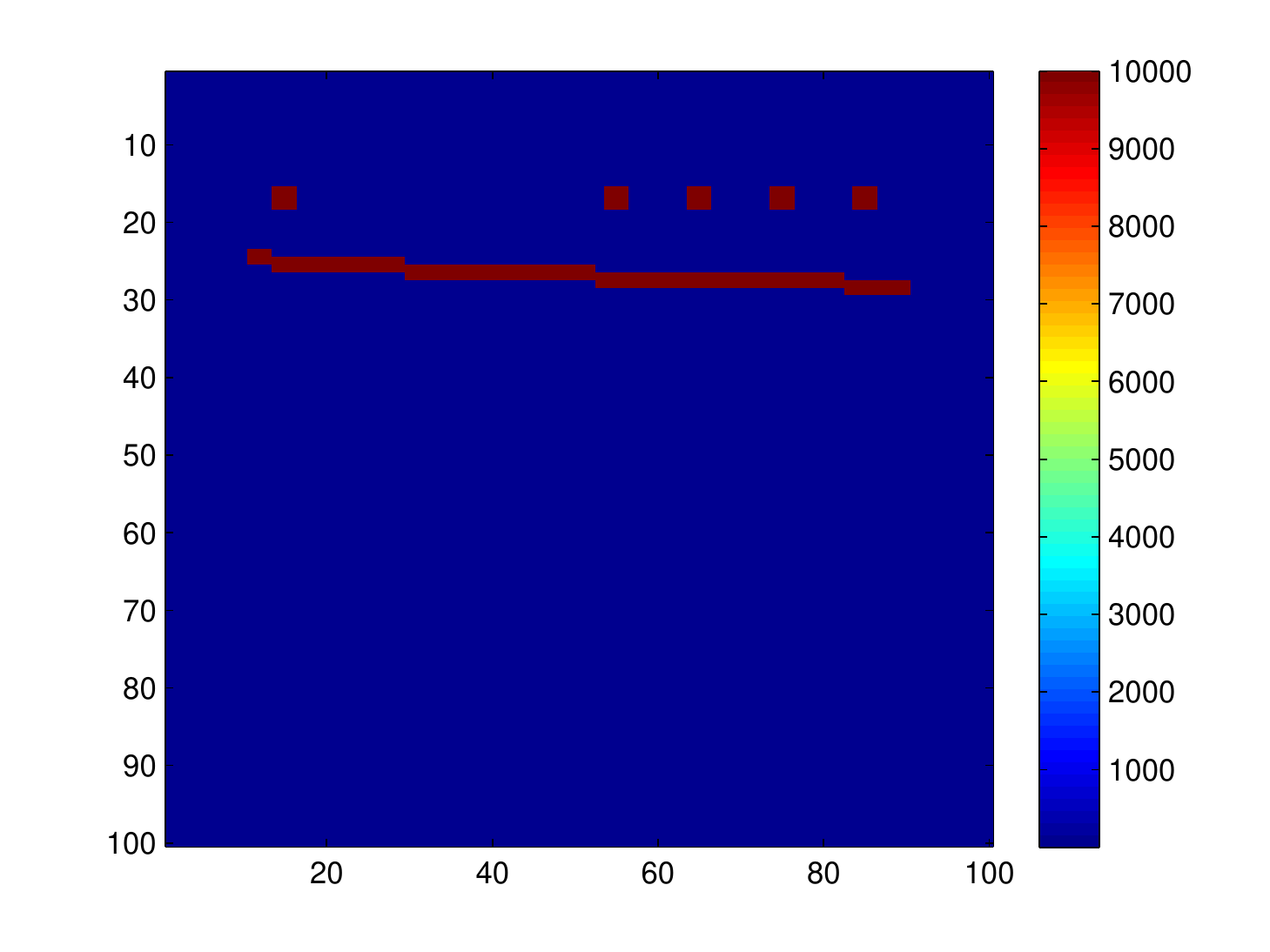}\includegraphics[scale=0.4]{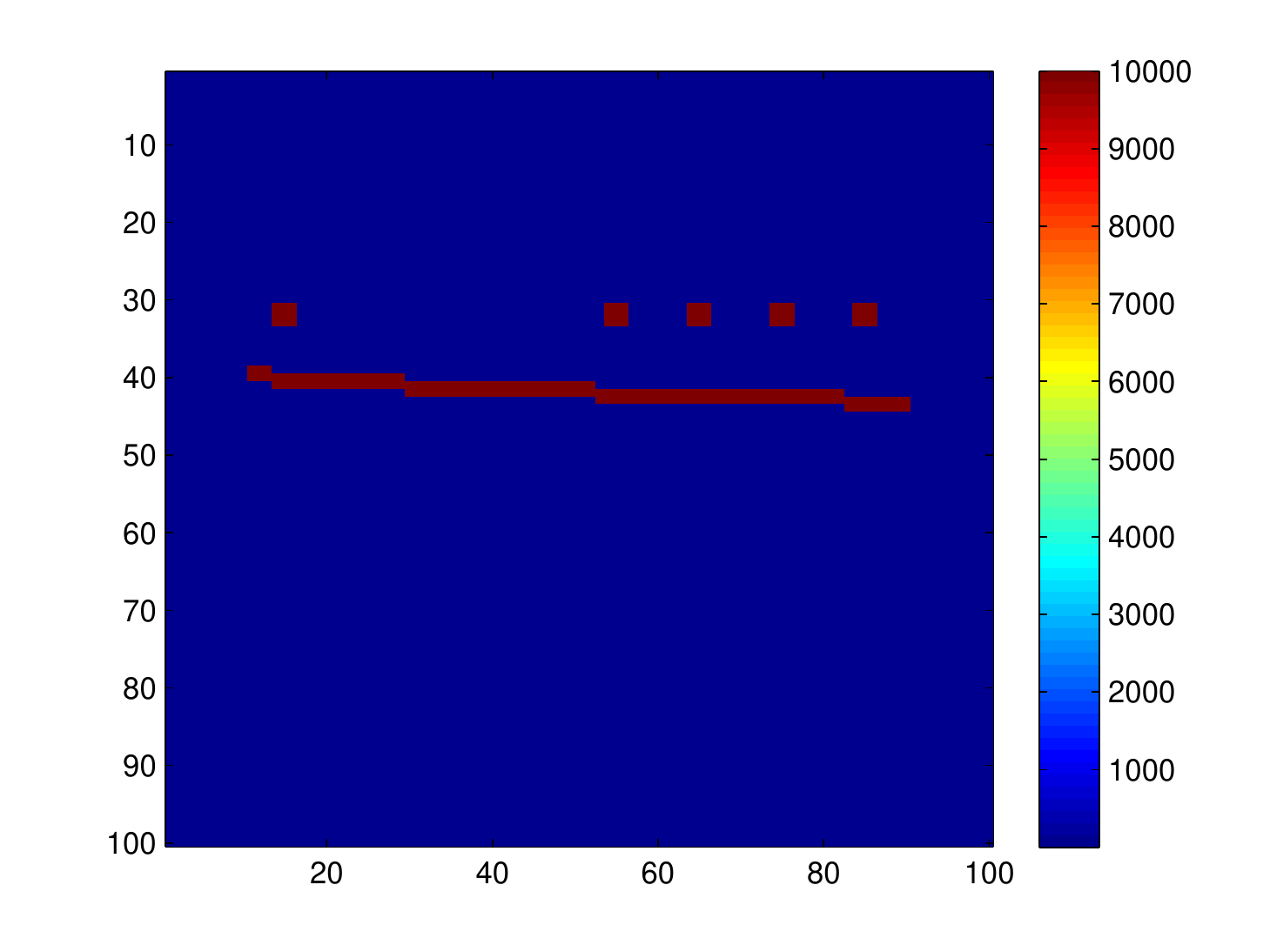}

\includegraphics[scale=0.4]{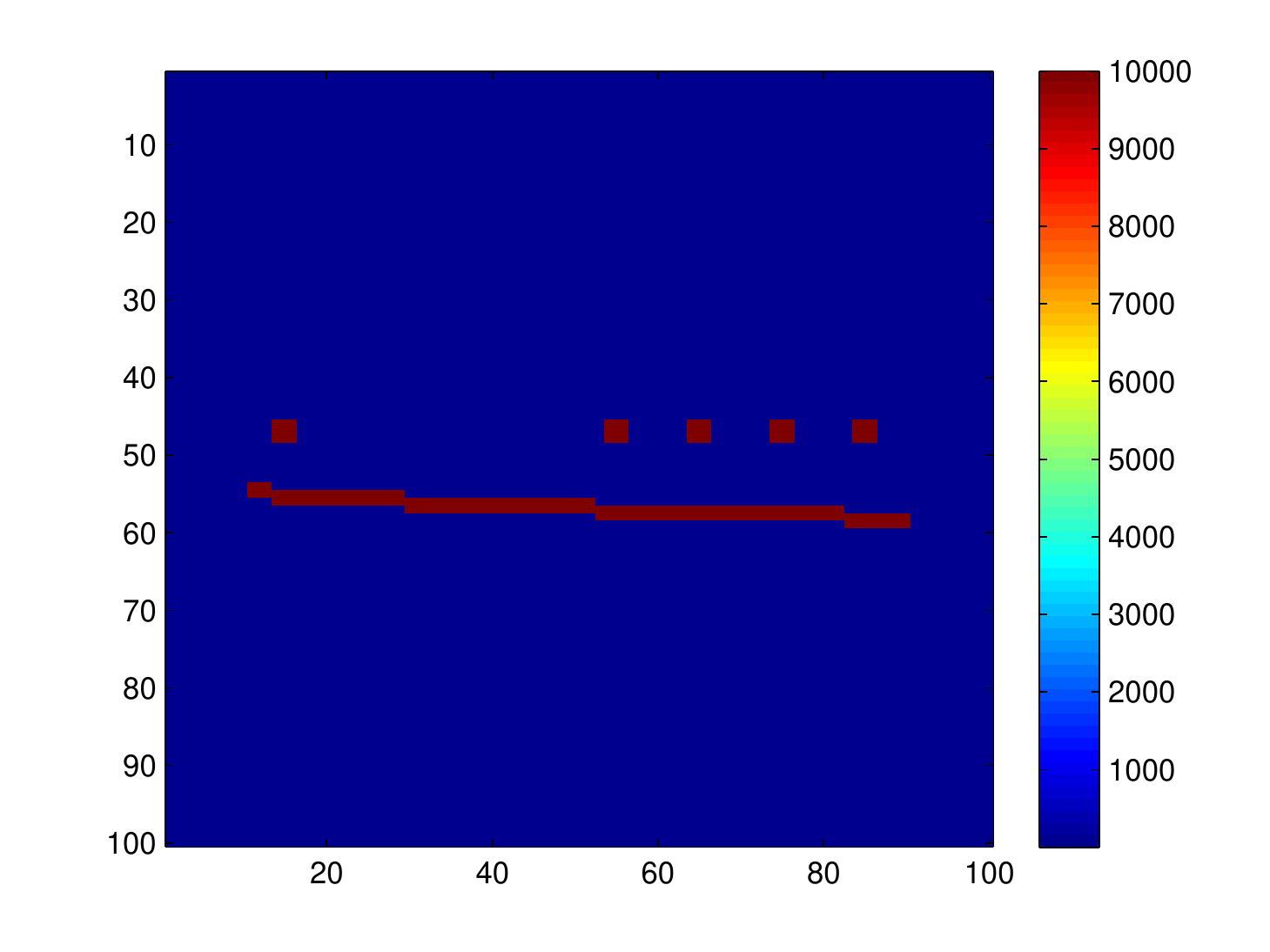}\includegraphics[scale=0.4]{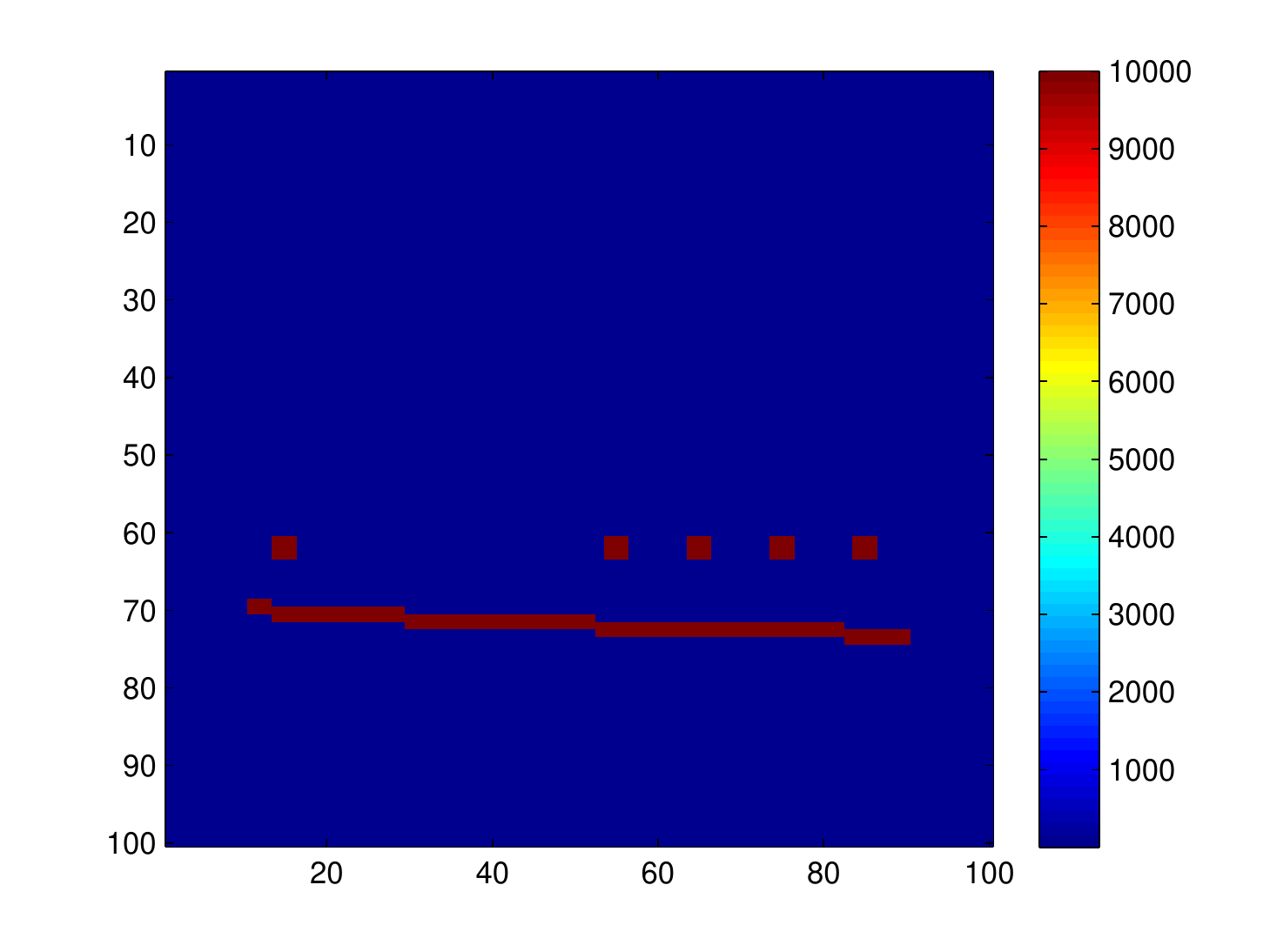}

\protect\caption{medium parameter. Top-Left: $\kappa(\mu_{1})$, Top-Right: $\kappa(\mu_{2}),$
Bottom-Left: $\kappa(\mu_{3})$, Bottom-Right: $\kappa(\mu_{4})$}
\label{fig:parameter_case2}
\end{figure}

{\bf Setup.} In our second example, we consider the source function $f$ to be the same as in the previous example. As for the medium parameter,  we use a nonlinear
function where the high permeability channel and inclusions move as
we change the parameter. The expression for the medium parameter is
\[
\kappa(\mu)= \kappa_1(x+\mu,y)
\]
In Figure \ref{fig:parameter_case2}, we show four different values of $\mu$, ($\mu=0,0.15,0.3,0.45$), that is used to construct the snapshot space. This is a complicated case as the
high-conductivity region is not fixed. 

\begin{figure}[ht]\centering
\includegraphics[scale=0.4]{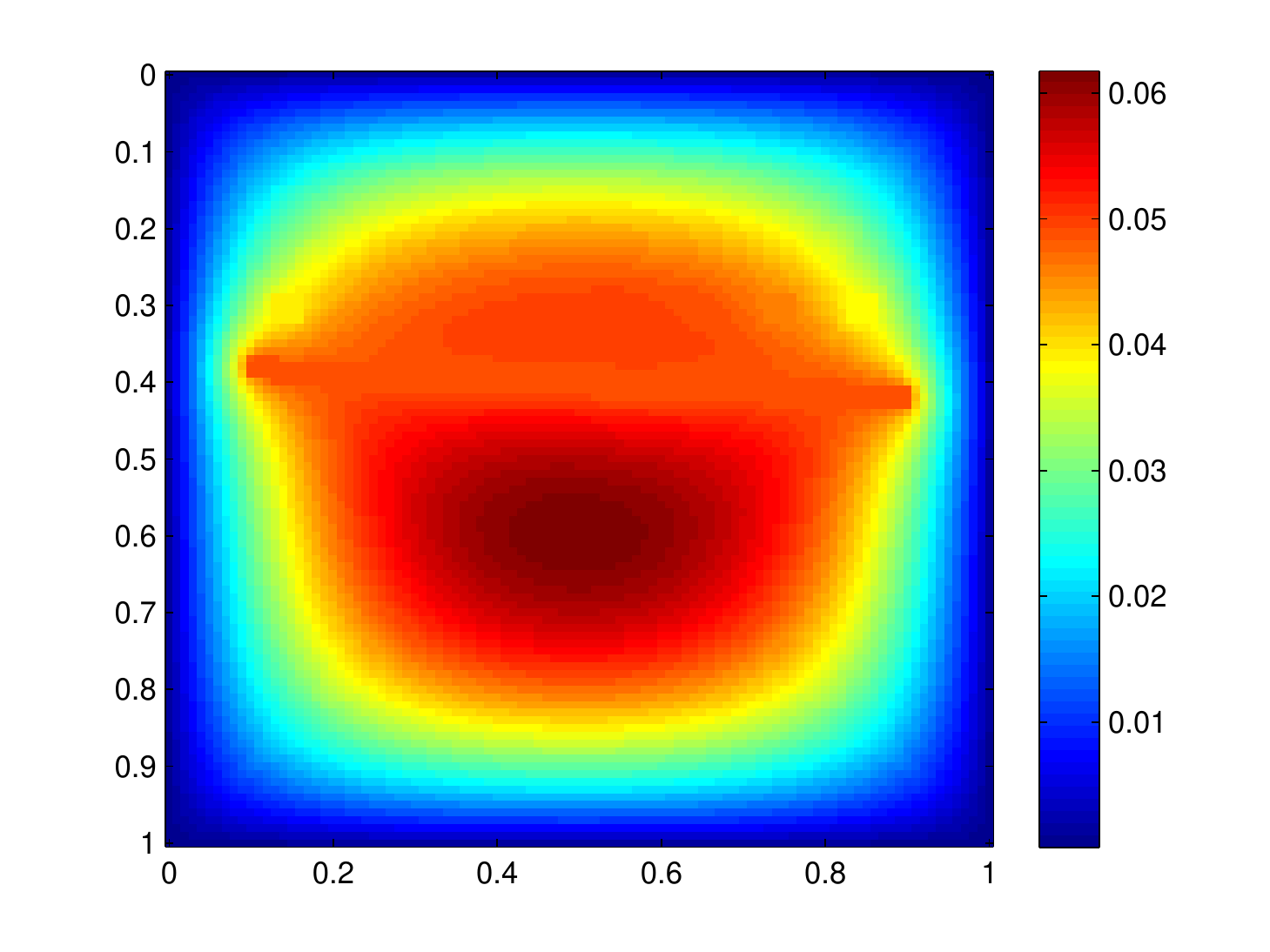}\includegraphics[scale=0.4]{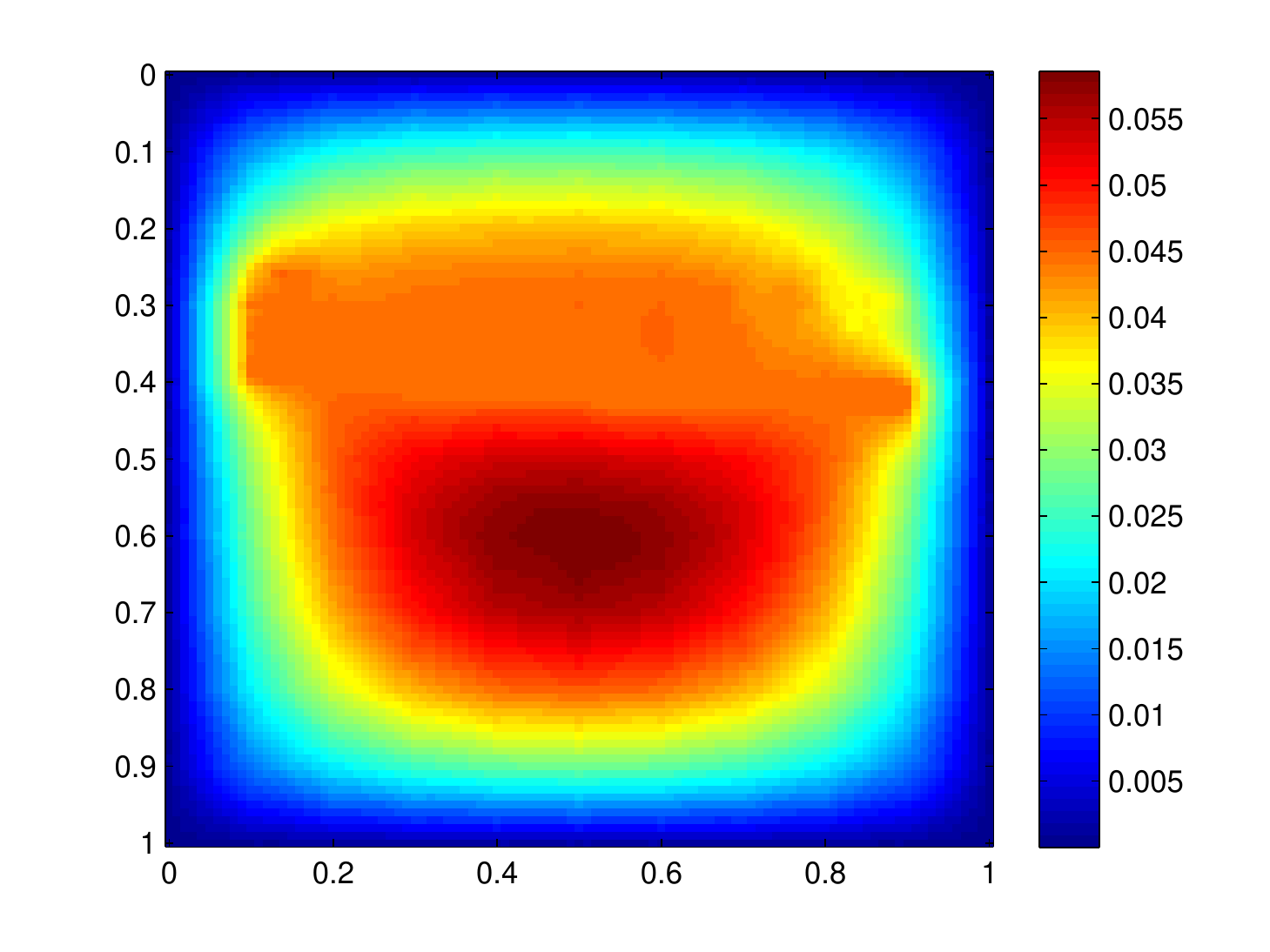}

\protect\caption{Left: Fine grid solution Right: Numerical solution ($14$ basis).}
\label{fig:sol_case2}
\end{figure}

{\bf Discussions of numerical results.} In Table \ref{table:sparse_DG Harmonic case2}, we show the convergence history of our method. The fine grid solution and the numerical solution are shown in Figure \ref{fig:sol_case2} with $\mu=0.14$. As we see from this table that the error is
larger compared to the previous case. The error decreases as we increase the dimension of the space. However, due to the fact that we do not span
all (or many) parameter values, the decay is slow.
As for the sparsity, we achieve a better sparsity compared to the previous
example because the online value of $\mu=0.14$ is close to one of
selected offline values $\mu=0.15$.
In fact, we observe
 that when the snapshot space dimension
is $9600$, as before, (i.e., $24$ randomized solutions per coarse block and per
each value of $\mu_i$), the nonzero coefficients in the expansion of
basis functions (over the whole domain) in terms of $9600$ snapshot
vectors are $3700$. The optimal expansion for this case (i.e.,
the case with $24$ randomized solutions per coarse block)
corresponds when $\mu=0.15$
is selected for snapshot construction and in this case, the number of
nonzero coefficients is $2400$ ($24$ per coarse region).
If we consider small dimensional online spaces and
the full snapshot space, then the sparsity is very small, as before.
For example, if we use $6$ randomized solutions per coarse block
and per each value of $\mu_i$ for identifying multiscale basis functions
per each coarse region, then, we will be using only $1/4$ of the snapshot
vectors and thus, the sparsity will be $9.5$ \%.
 As we observe
that our numerical examples identify appropriate sparsity of the solution
space. Again,
we expect a more significant gain in the sparsity if more parameter
values are used.

\begin{table}[htb!]
\centering

\begin{tabular}{|c|c|c|c|c|}
\hline
\multirow{2}{*}{$\text{dim}(V_{\text{on}})$}  &
\multicolumn{2}{c|}{  $\|u-u^{\text{on}} \|$ (\%) }  \\
\cline{2-3} {}&
$\hspace*{0.8cm}   L^{2}(D)   \hspace*{0.8cm}$ &
$\hspace*{0.8cm}   H^{1}_\kappa(D)  \hspace*{0.8cm}$
\\
\hline\hline
      $800$    &    $13.50$    & $31.06$ \\
\hline
     $1000$   &   $12.02$       & $29.45$  \\
\hline
     $1200$     &   $9.72$    & $26.66$  \\
\hline
      $1400$    &    $7.97$    & $24.13$ \\
\hline
     $\text{dim}(V_{\text{snap}})$   &   $5.78$       & $20.27$  \\
\hline
\end{tabular}
\caption{Convergence history of the DGMsFEM using oversampling Harmonic basis. The fine-scale dimension is 12100. The full snapshot space dimension is 9600. }
 \label{table:sparse_DG Harmonic case2}
\end{table}

\begin{remark}
We have implemented CG-GMsFEM using spectral basis approach and observed similar results.

\end{remark}

\subsection{Second Approach. Sparse Snapshot Subspace Approach}

In this section, we will show the numerical example by using second approach to directly calculate the sparse multiscale solution by $l_1$ minimization.

\subsubsection*{Example 1}

\begin{figure}[ht]
\centering

\includegraphics[scale=0.4]{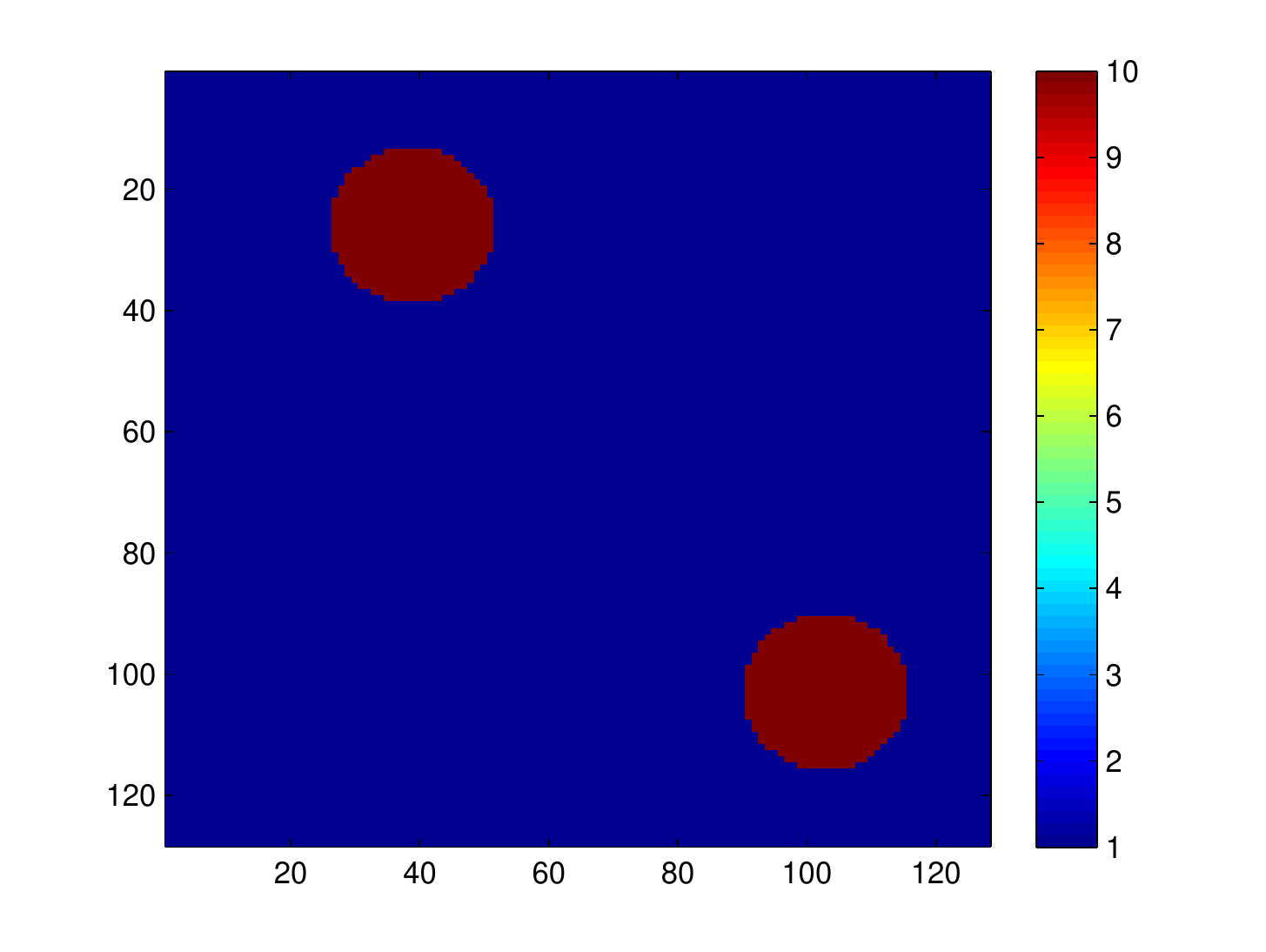}

\protect\caption{Parameter $n(x)$.}
\label{fig:para_n}
\end{figure}

{\bf Setup.} In this example, we consider the domain $D=[0,1]^{2}$
that is partitioned into the
coarse grid with grid size $H=1/8$ and each coarse block is subdivided
into $16\times16$ fine square blocks with length $h=H/16$.
Therefore, the fine mesh size $h=\cfrac{1}{128}$. We consider $\Omega=2$, $\kappa\equiv1$ and $n(x)$
is shown in Figure \ref{fig:para_n}. We consider a zero source function with Dirichlet boundary condition $g$, given by $g=e^{-i\Omega k\cdot x}$ where $k=(\sin(\cfrac{\pi}{4}),\cos(\cfrac{\pi}{4}))$.

\begin{figure}[htb!]
\centering

\includegraphics[scale=0.4]{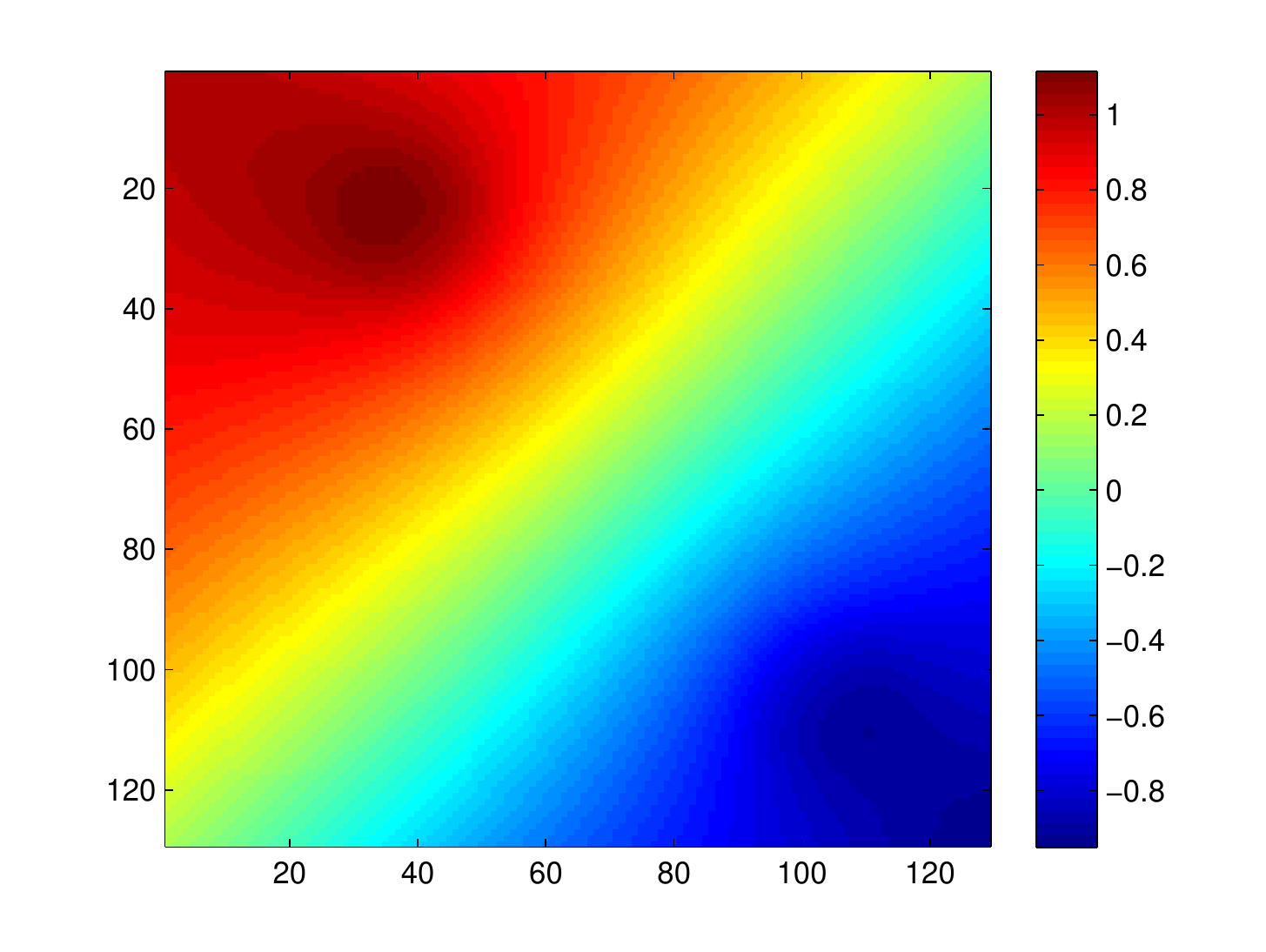} \includegraphics[scale=0.4]{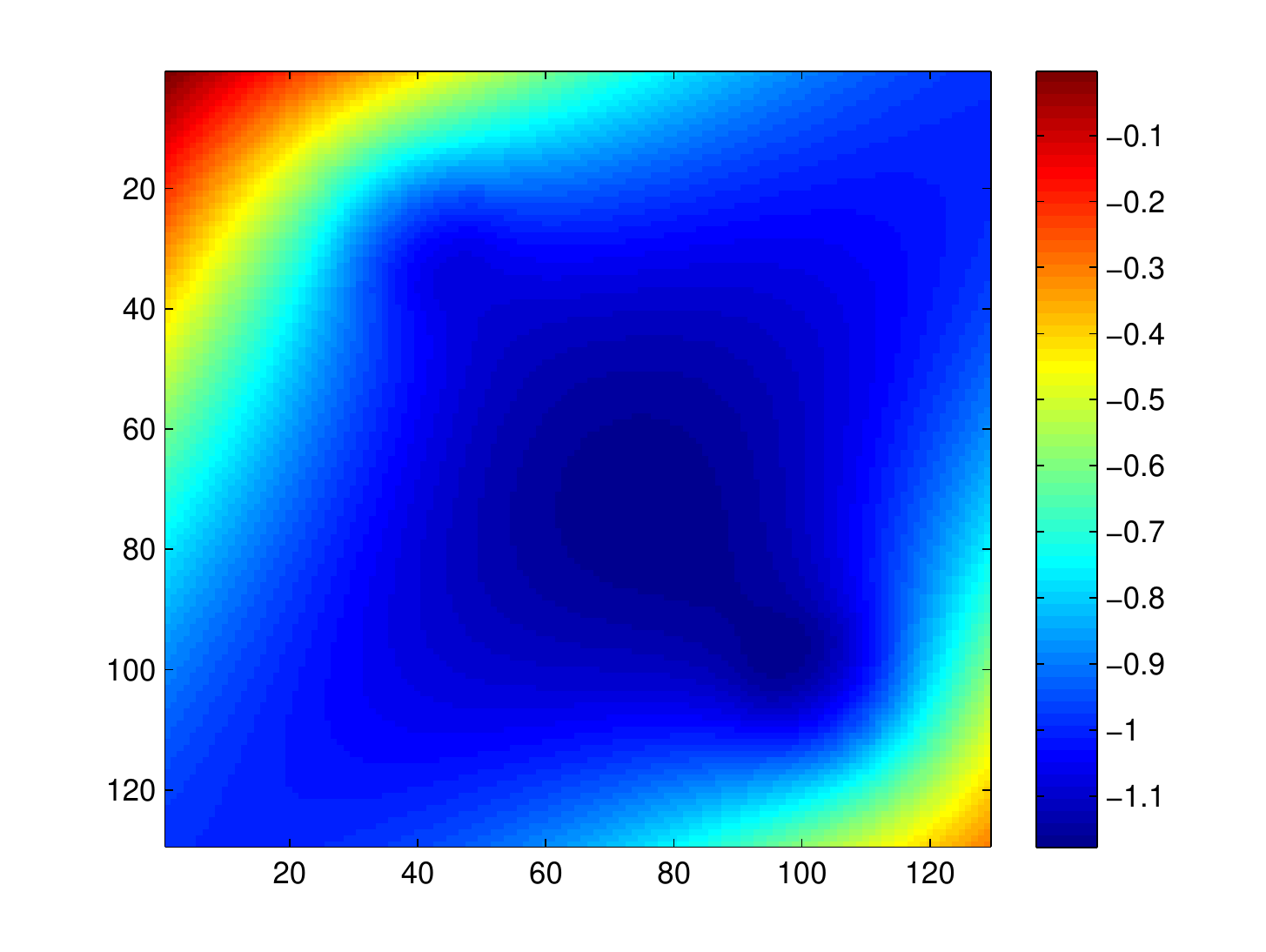}

\protect\caption{Reference solution $u$, Left: Real part of the solution; Right: Imaginary part of the solution.}

\label{fig:case2_sol}
\end{figure}

\begin{figure}[ht]
\centering

\includegraphics[scale=0.4]{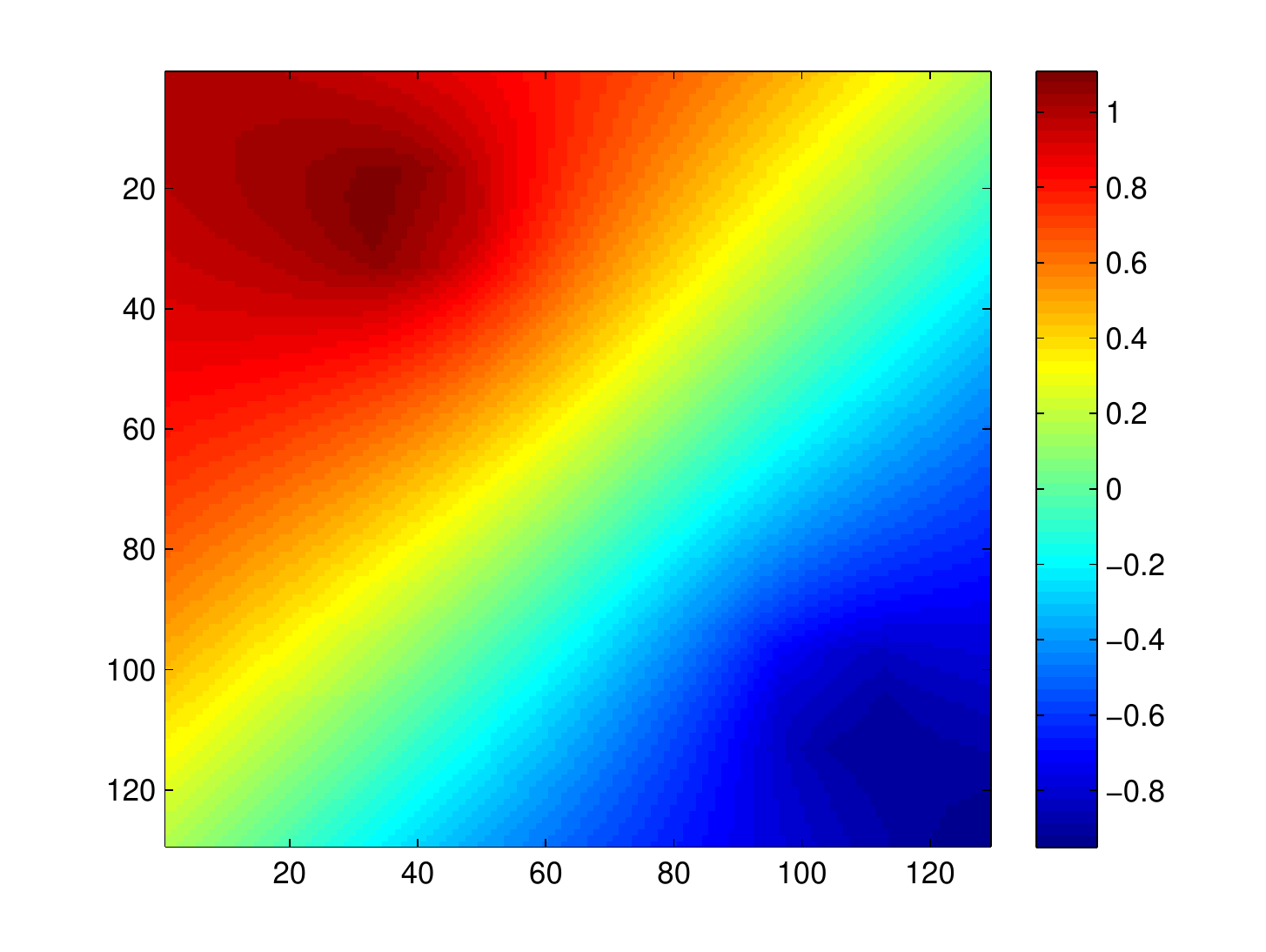} \includegraphics[scale=0.4]{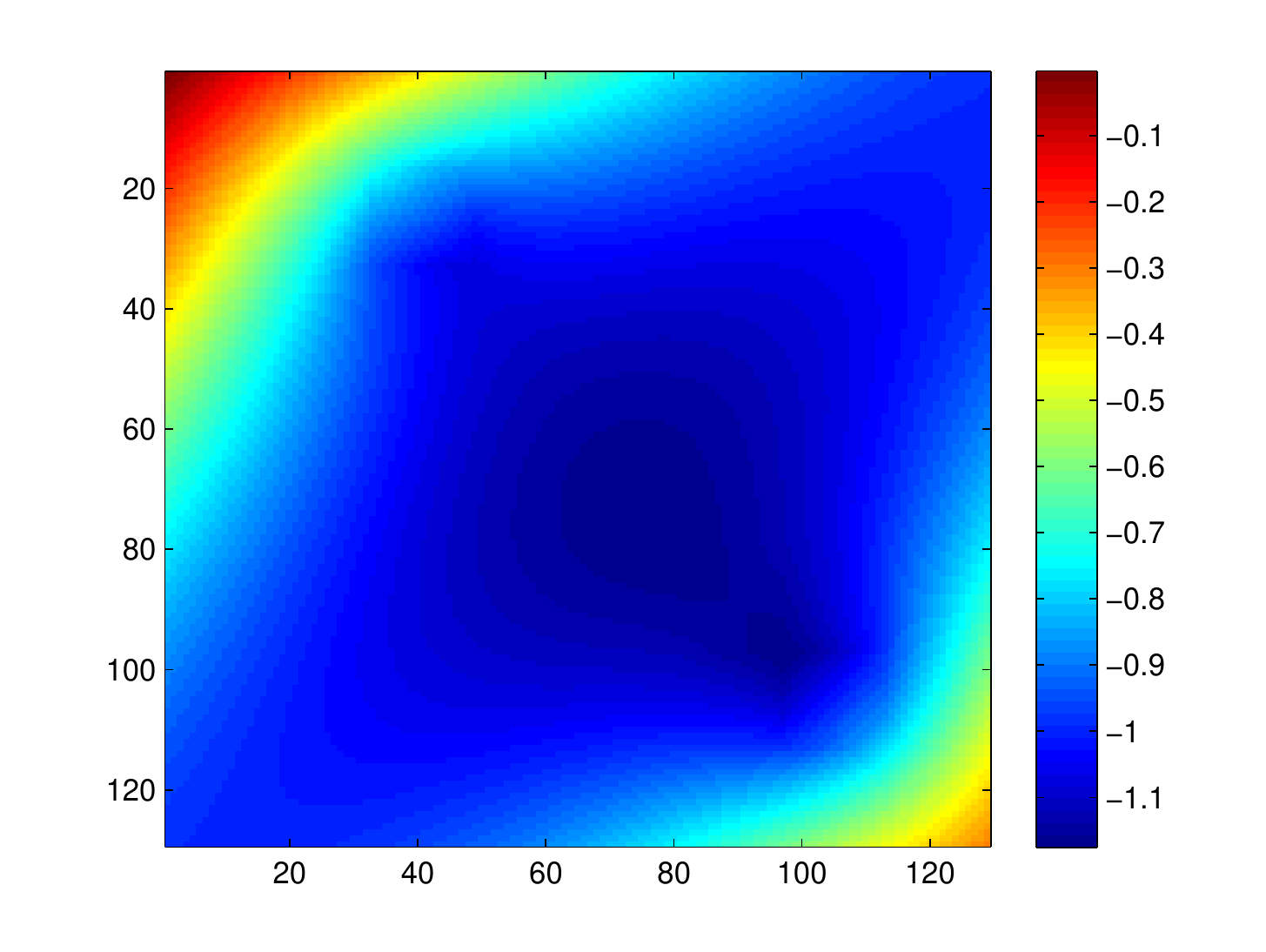}

\protect\caption{Snapshot solution (i.e., when using all snapshot vectors), $u$, Left: Real part; Right: Imaginary part.}
\label{fig:case2_snap}
\end{figure}

\begin{figure}[ht]
\centering

\includegraphics[scale=0.4]{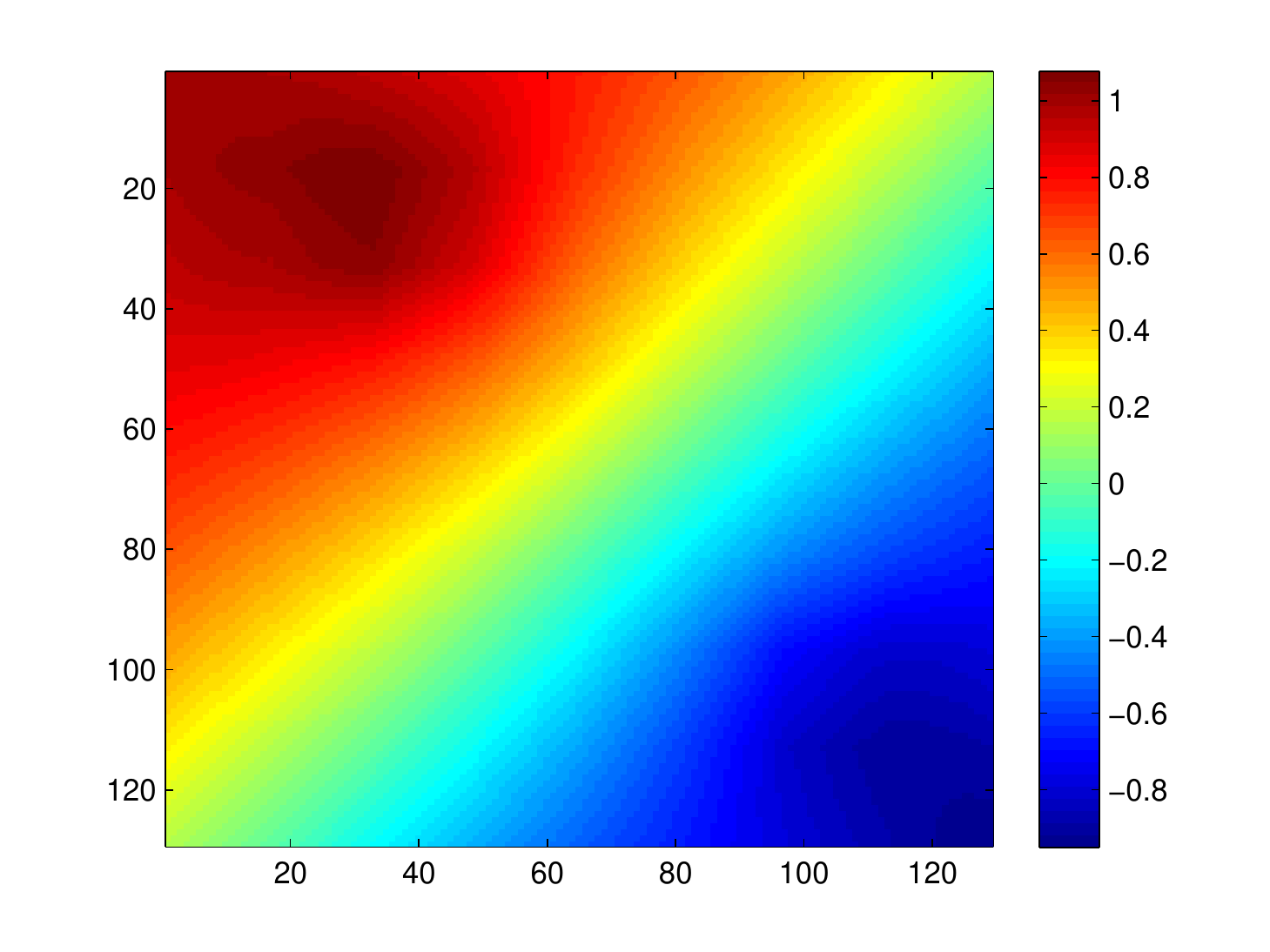} \includegraphics[scale=0.4]{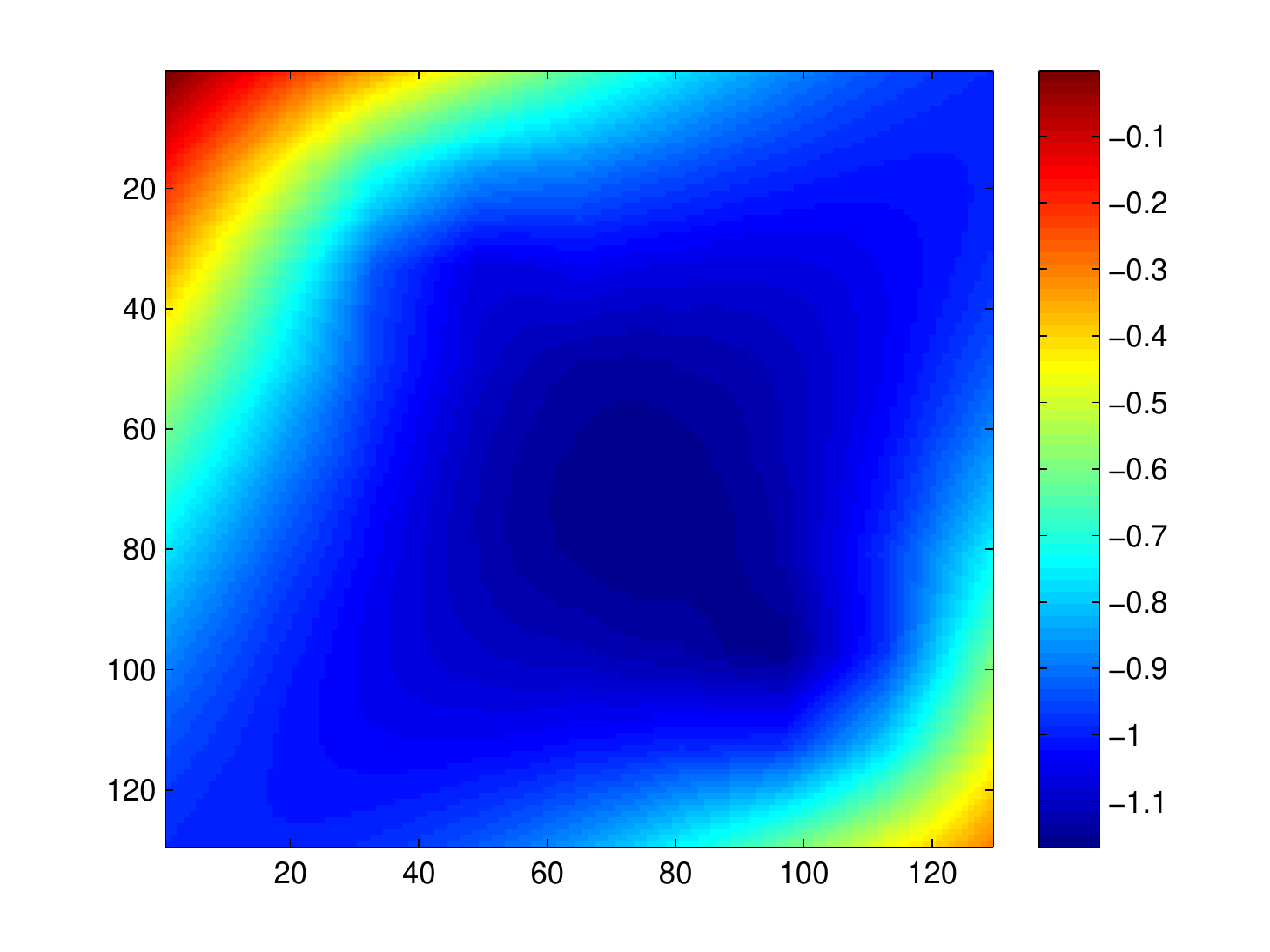}

\protect\caption{Numerical solution with $4$ test basis functions, $u$, Left: Real part of the solution. Right: Imaginary part of the solution.}
\label{fig:case2_test4}
\end{figure}

{\bf Discussions of numerical results.}
We will compare our result with the reference solution, which is calculated on the fine grid and shown in Figure
\ref{fig:case2_sol}.
Notice that, within each coarse grid block, the reference solution has few dominant propagating directions,
which suggests sparsity of the solution in the snapshot space.
In this case, the snapshot space is spanned by  local plane waves with dimension $\text{dim}(V_\text{snap})=1280$,
as defined in (\ref{eq:SparseTrig}) with $k_i$'s distributed uniformly.
 The snapshot solution (i.e., if we use all snapshot vectors)
has
 $1.63\%$ relative error with respect to the fine-scale solution.
We compare the solutions in Figure \ref{fig:case2_snap}.
As we observe, the snapshot solution is accurate.
Next, we calculate the sparse solution by varying the dimension of the test
space.
The latter defines a sparse solution in the subspace of the test space.
The numerical solution calculated with $4$ test basis per coarse grid block is shown in Figure \ref{fig:case2_test4}. In Table \ref{table:sparse_DG case2}, we show the convergence history of the second approach, where $\|u\|_{H^1(D)}^2 = \int_D |\nabla u|^2$.
As we observe that for low dimensional test spaces, the solution is very
sparse in the snapshot space
(and this sparsity is
about the same as the test space). We increase the dimension of the test space
to achieve a higher accuracy.
For example, for the solution
with the sparsity $408$ (i.e., $408$ non-zero coefficients in the span of $1280$
snapshot vectors), we have $1.63$ \% $L^2$ error.

{\bf Why to expect a sparsity.} Next, we briefly describe why to expect
a sparsity in this problem.
The snapshot space consists of
local problems corresponding to different directions $k_m$ (see (\ref{eq:SparseTrig})). In this problem,
we expect that the solution will consist of plane wave solutions with a few
directions. By using plane wave solutions, we can identify these few directions. Note that in this example, we can not identify local spectral decomposition.


\begin{table}[htb!]
\centering

\begin{tabular}{|c|c|c|c|c|c|}
\hline
\multirow{2}{*}{$\text{dim}(V_{\text{test}})\times 2$}  &
\multicolumn{2}{c|}{  $\|u-u_{\text{ms}} \|$ (\%) }  &\\
\cline{2-3} {}&
$\hspace*{0.8cm}   L^{2}(D)   \hspace*{0.8cm}$ &
$\hspace*{0.8cm}   H^{1}(D)  \hspace*{0.8cm}$ &
sparsity of the sol
\\
\hline\hline
       $128$     &   $41.10$    & $195.38$ &128 \\
\hline
      $256$    &    $21.12$    & $45.16$  &252\\
\hline
      $384$    &    $14.62$    & $32.62$  &344\\
\hline
     $512$   &   $3.49$       & $14.40$   &408\\
\hline
     $\text{dim}(V_{\text{snap}})$   &   $1.63$       & $9.88$ & 1280\\
\hline
\end{tabular}
\caption{Convergence history of the DGMsFEM to compute sparse multiscale solution directly. The fine-scale dimension is 16384. The snapshot space dimension is 1280.}
 \label{table:sparse_DG case2}
\end{table}

\subsubsection*{Example 2}

\begin{figure}[ht]
\centering
\includegraphics[scale=0.4]{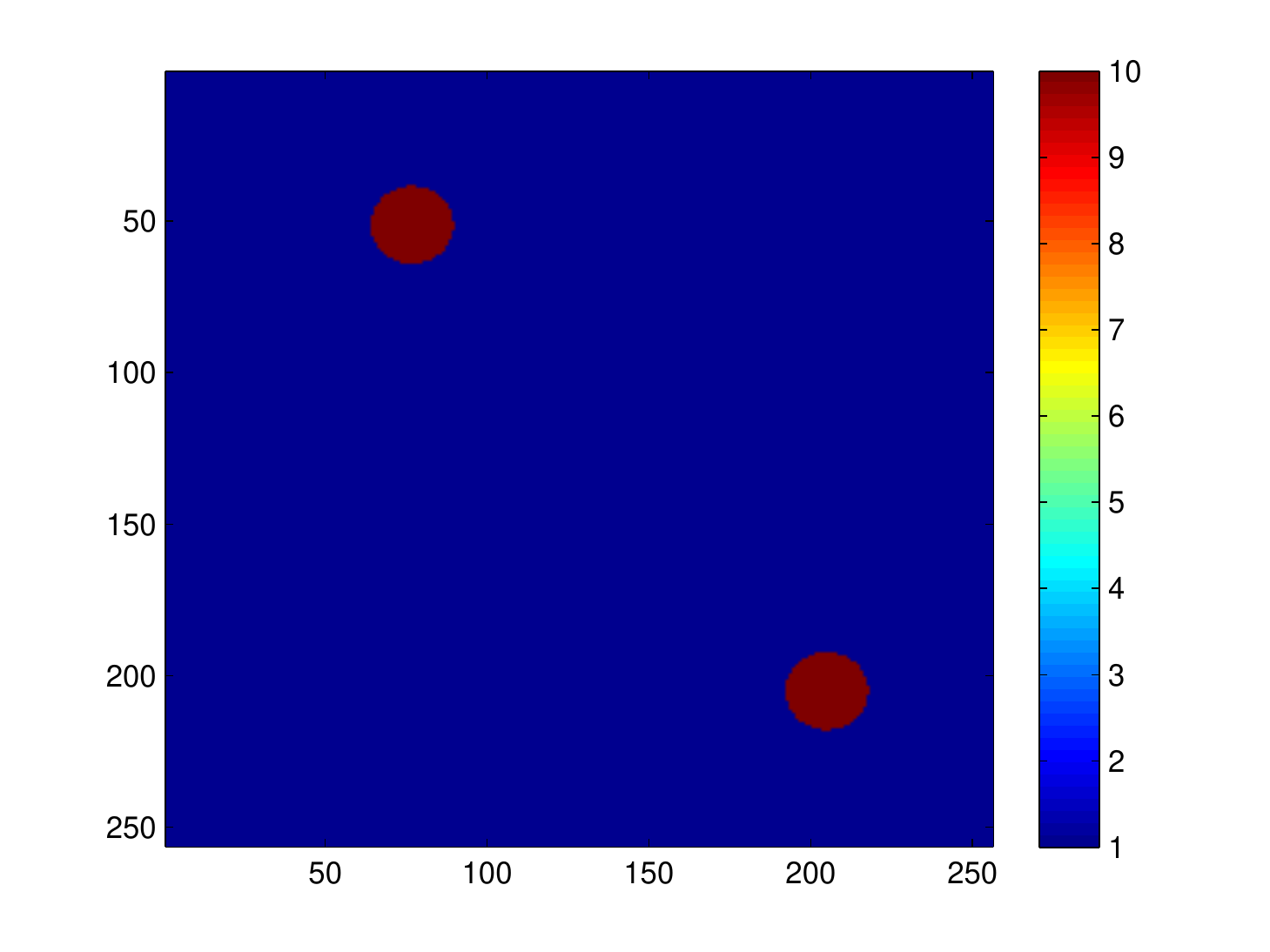}
\protect\caption{Parameter $n(x)$.}
\label{fig:para_n2}
\end{figure}

{\bf Setup.} In this example,
we consider the domain $D=[0,1]^{2}$ is partitioned into the
coarse grid with grid size $H=1/16$ and each coarse block is subdivided
into $16\times16$ fine square fine block with length $h=H/16$, therefore, the fine mesh size $h=\cfrac{1}{128}$. We consider $\Omega=8$ and $n(x)$
is shown in figure \ref{fig:para_n2}. The parameter $\kappa$, source function $f$, and boundary condition $g$, are the same as the previous example. Because of higher value of $\Omega$,
we take the fine grid $2$ times finer.

\begin{figure}[htb!]
\centering

\includegraphics[scale=0.4]{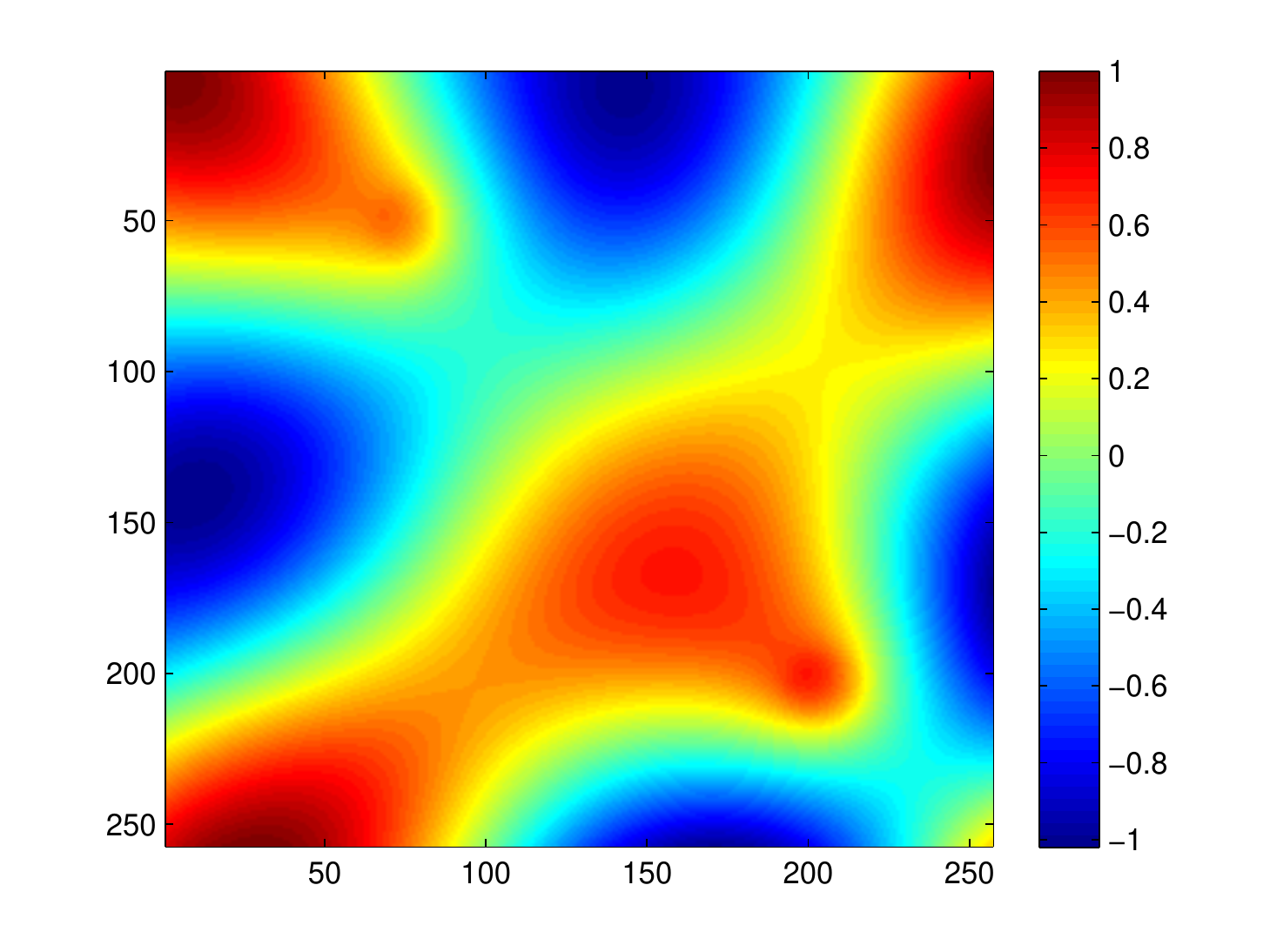} \includegraphics[scale=0.4]{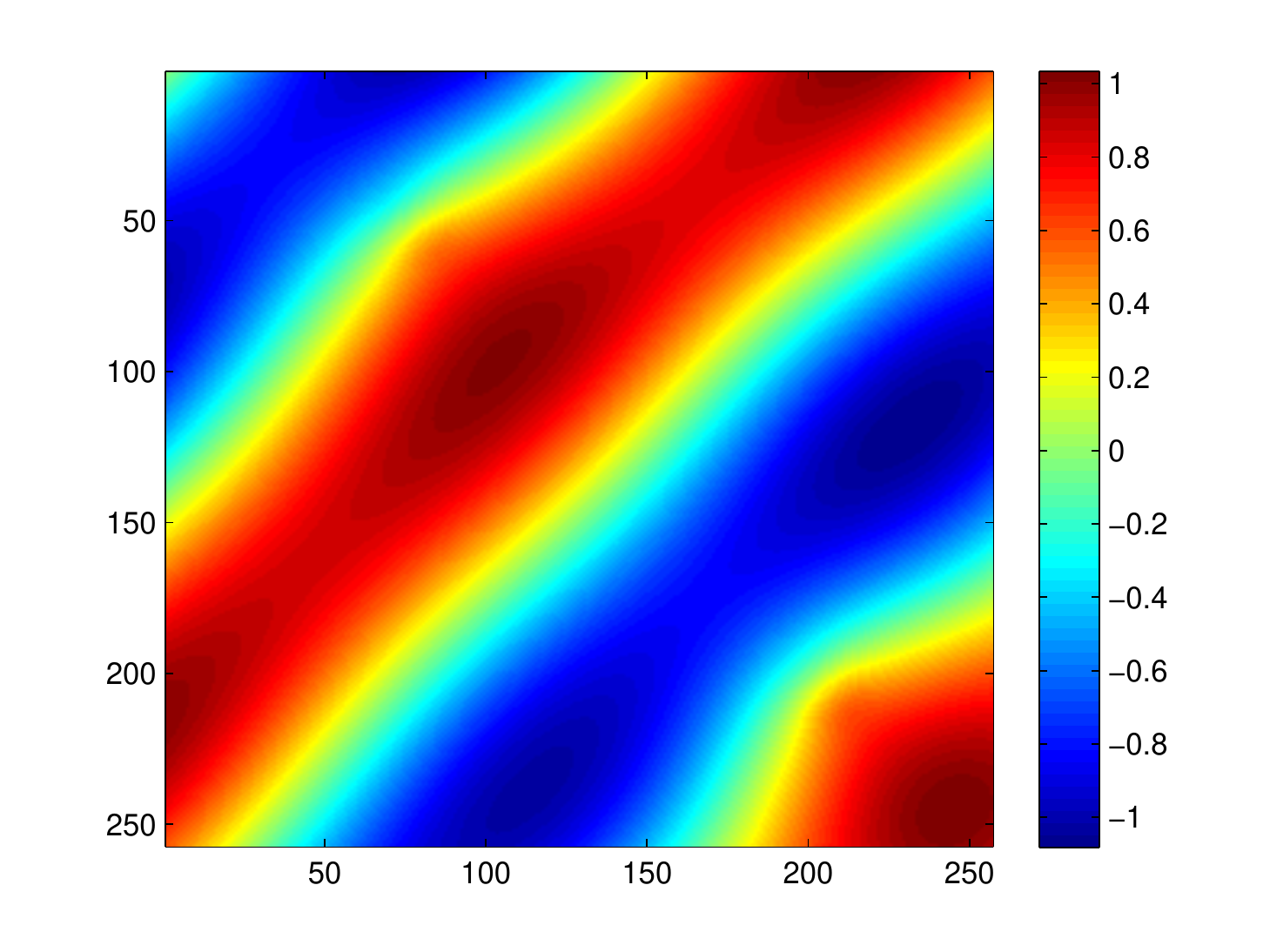}

\protect\caption{Reference solution $u$, Left: Real part of the solution; Right: Imaginary part of the solution.}

\label{fig:case2_sol2}
\end{figure}

\begin{figure}[ht]
\centering

\includegraphics[scale=0.4]{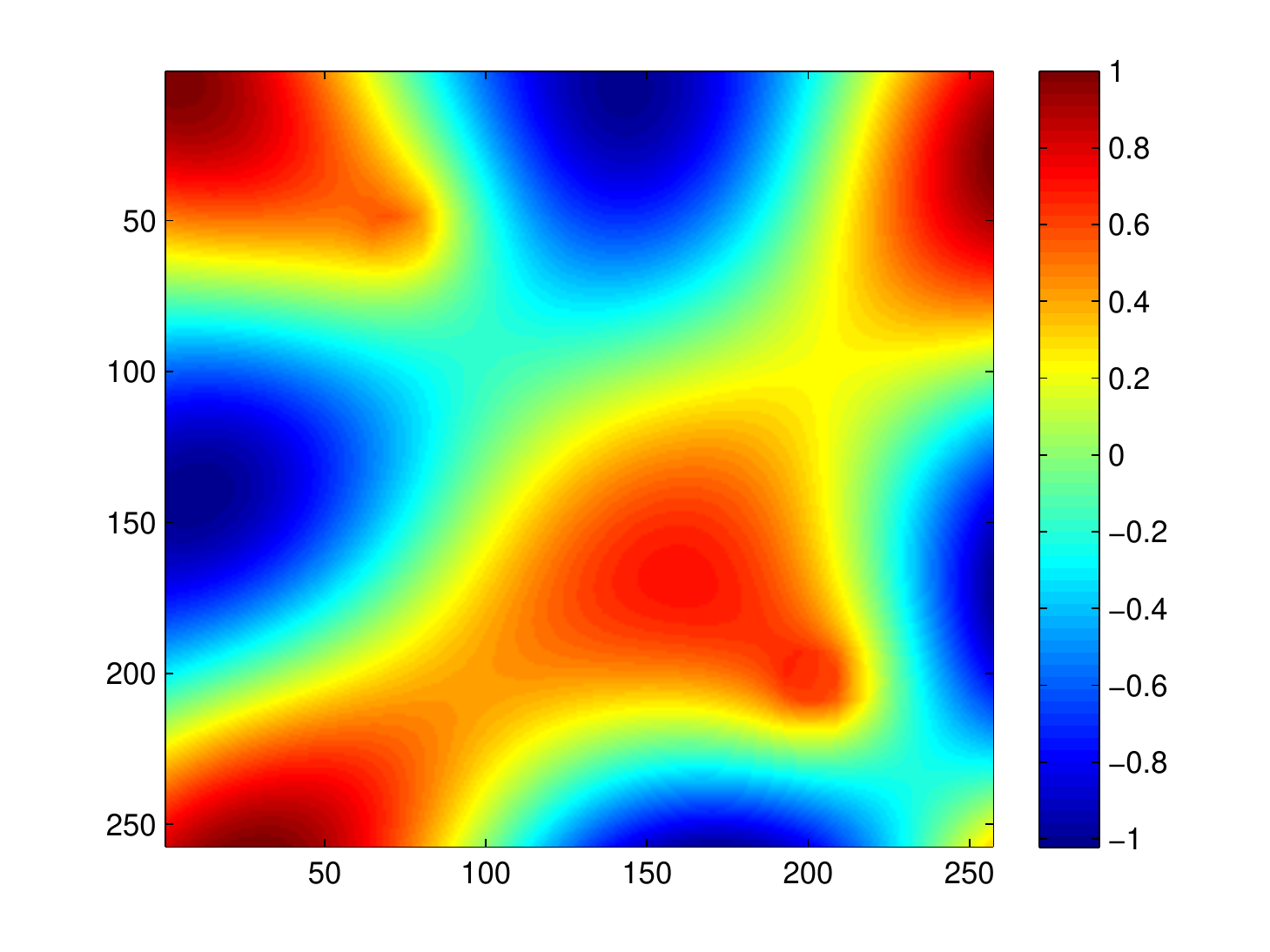} \includegraphics[scale=0.4]{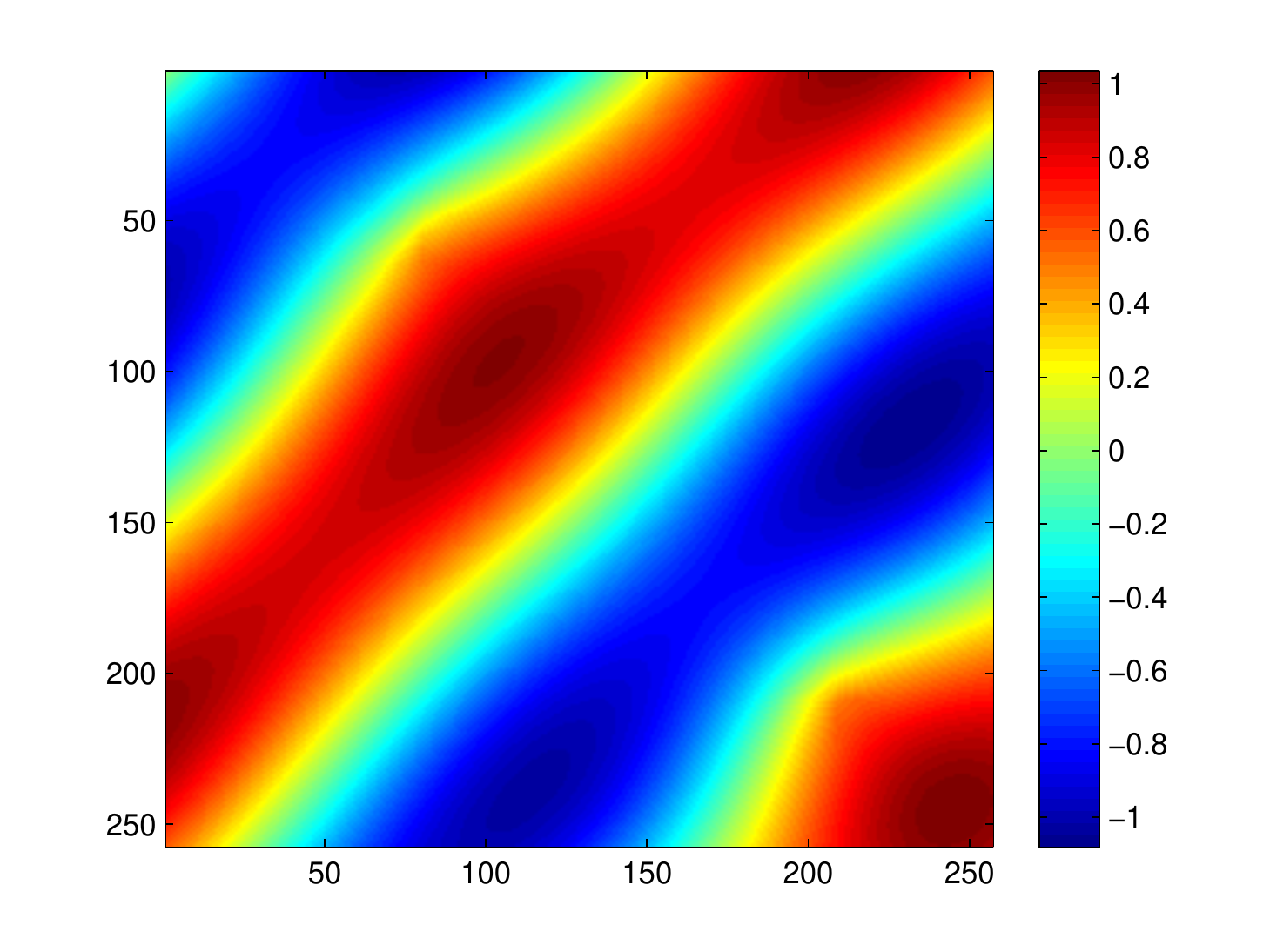}

\protect\caption{Snapshot solution (i.e., when using all snapshot vectors), $u$, Left: Real part; Right: Imaginary part.}
\label{fig:case2_snap2}
\end{figure}

\begin{figure}[ht]
\centering

\includegraphics[scale=0.4]{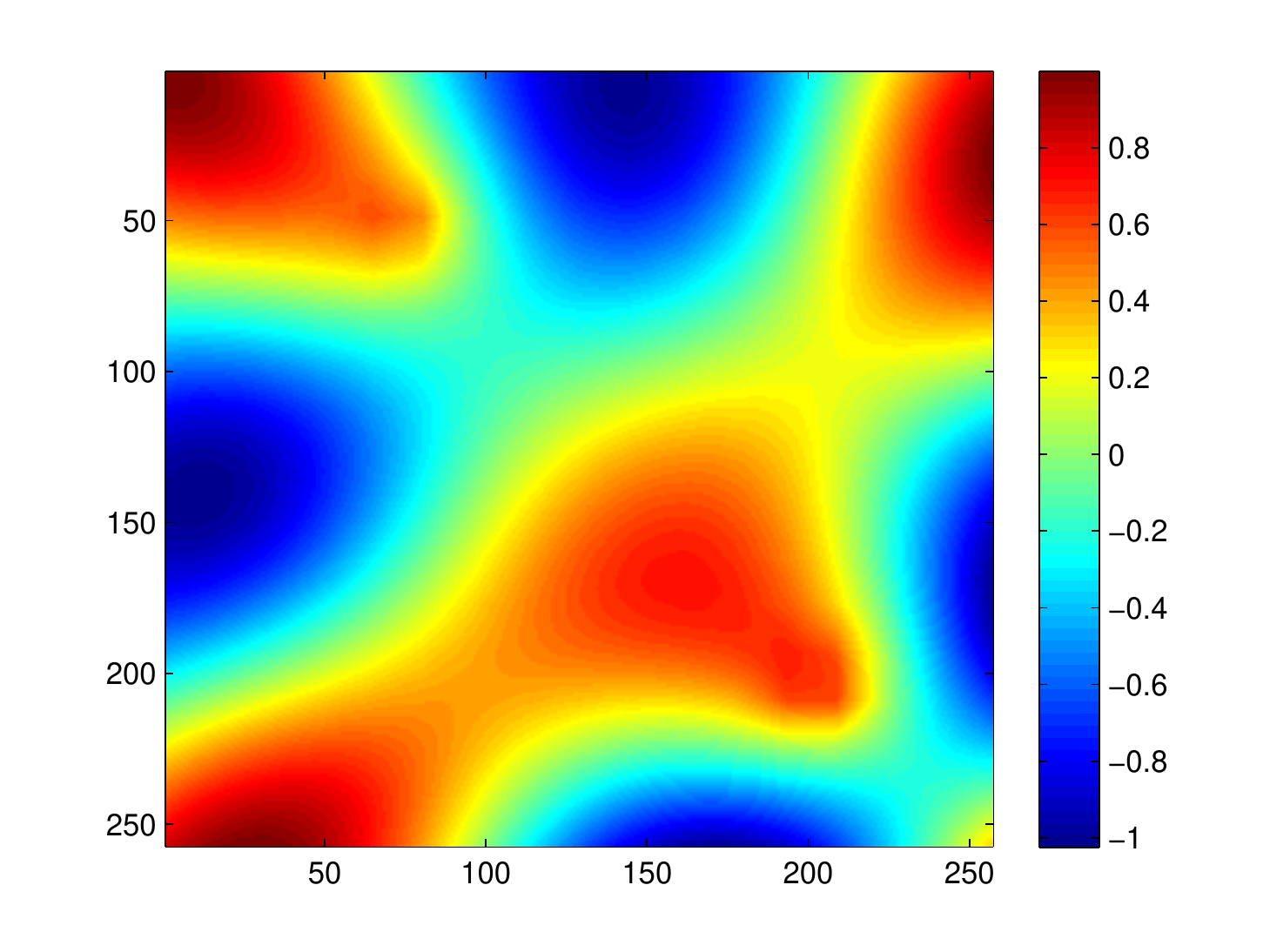} \includegraphics[scale=0.4]{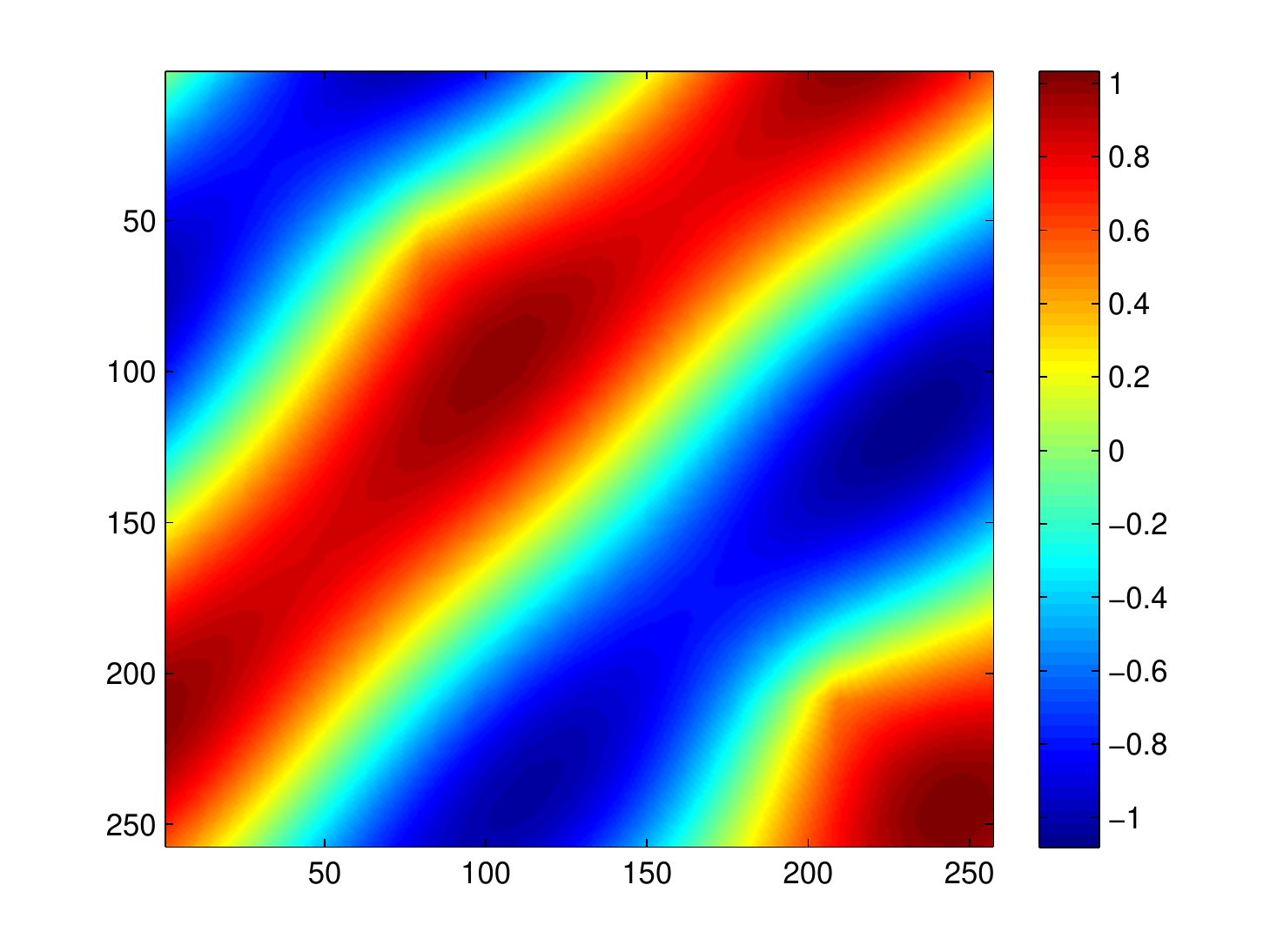}

\protect\caption{Numerical solution with $4$ test basis functions, $u$, Left: Real part of the solution. Right: Imaginary part of the solution.}
\label{fig:case2_test5_2}
\end{figure}

{\bf Discussions of numerical results.}
We will compare our results with the multiscale approach with the reference solution, which calculated on the fine grid and shown in Figure \ref{fig:case2_sol2}.
 Notice that, within each coarse grid block, the reference solution has few dominant propagating directions,
which suggests sparsity of the solution in the snapshot space.
 In this case, the snapshot space is spanned by local plane waves with dimension $\text{dim}(V_\text{snap})=5120$,
as defined in (\ref{eq:SparseTrig}) with $k_i$'s distributed uniformly.
The snapshot solution error is  $2.44\%$ and it is shown in Figure \ref{fig:case2_snap2}.
As we observe the snapshot solution is accurate.
Next, we calculate the sparse solution by varying the dimension of the test
space.
The latter defines a sparse solution in the subspace of the test space.
The numerical solution calculated with $4$ test basis per coarse grid block is shown in Figure \ref{fig:case2_test5_2}.
In Table \ref{table:sparse_DG case2_2}, we show the convergence history of the second approach.
As we observe that for low dimensional test spaces, the solution is very sparse
in the snapshot space. We increase the dimension of the test space to
achieve a higher accuracy.
 For example, the solution with $1958$ nonzero coefficients in the snapshot
space provides $4.25$ \% accurate solution in $L^2$ sense.



\begin{table}[htb!]
\centering

\begin{tabular}{|c|c|c|c|c|c|}
\hline
\multirow{2}{*}{$\text{dim}(V_{\text{test}})\times 2$}  &
\multicolumn{2}{c|}{  $\|u-u_{\text{ms}} \|$ (\%) }  &\\
\cline{2-3} {}&
$\hspace*{0.8cm}   L^{2}(D)   \hspace*{0.8cm}$ &
$\hspace*{0.8cm}   H^{1}(D)  \hspace*{0.8cm}$ &
sparsity of the sol
\\
\hline\hline
       $512$     &   $79.12$    & $176.91$ &500 \\
\hline
      $1024$    &    $74.07$    & $103.93$  &913\\
\hline
      $1536$    &    $40.69$    & $51.79$  &1294\\
\hline
     $2048$   &   $17.27$       & $23.74$   &1591\\
\hline
     $2560$   &   $4.25$       & $8.63$   &1958\\
\hline
     $\text{dim}(V_{\text{snap}})$   &   $2.44$       & $6.18$ & 5120\\
\hline
\end{tabular}
\caption{Convergence history of the DGMsFEM to compute sparse multiscale solution directly. The fine-scale dimension is 66049. The snapshot space dimension is 5120.}
 \label{table:sparse_DG case2_2}
\end{table}

\section{Conclusions}
\label{sec:conclusion}

\subsection{Summary of the results}

In the paper, we develop approaches to identify sparse multiscale basis functions
in the snapshot space within GMsFEM.
The snapshot spaces are constructed in a special way that
allows sparsity for the solution. We consider two apporaches. In the first approach,
local multiscale basis functions are constructed, which are sparse in the snapshot
space. These multiscale basis functions are constructed by identifying dominant
modes in the snapshot space using $l_1$ minimization techniques.
As for the application, we consider parameter-dependent multiscale problems.
In the second approach, we apply $l_1$ minimization techniques directly to solve
the global problem. This approach is more expensive as it directly deals with
a large snapshot space. As for the application, we consider Helmholtz equations.
For both approaches and their respective applications, we present numerical
results and discuss computational savings. Our numerical examples are
simplistic and are designed to convey the main idea of the proposed approach.

\subsection{Sparsity assumption}

Both approaches assume that the solution is sparse in the snapshot space.
The latter requires special snapshot spaces, which can yield this sparsity.
For example. for local snapshot vectors considered in the paper, this requires
identifying boundary conditions for the snapshot solutions, which can sparsily
represent
the solution. This may not be easy in general, though in some examples can still
be achieved. Besides examples presented in the paper, one can consider
scale separation cases and use piecewise linear boundary conditions.
Whether a general framework for constructing these snapshot vectors is possible
remains an open question.

\subsection{Adaptivity and online basis functions}

In general, once offline spaces are identified,
we can use adaptivity \cite{chung2014adaptive1, chung2014adaptive2}
and online basis functions
\cite{chung2015residual,chung2015online} to achieve a small error.
The adaptivity is accomplished by identifying the regions
with large residuals and enriching the spaces in those regions.
In our earlier works \cite{chung2014adaptive1, chung2014adaptive2},
we have shown that one needs to use some special
error indicators. For the first approach, we can use the ``next''
eigenvector obtained from local eigenvalue decomposition
to construct multiscale basis functions. For the second
approach, one can increase the test space for additional
multiscale basis functions.

In the regions with largest residuals,
we can also use online basis functions to reduce the error substantially.
Online basis functions are computed locally and identified as
the localized basis functions which can give a largest reduction
in the error. These basis functions involve solving local
problems with a residual on the right hand side (see \cite{chung2015online} for
online basis functions for DG).

In this paper, we can apply adaptivity and online basis functions
as discussed in \cite{chung2014adaptive2, chung2015online}.
For parameter-dependent problems,
one can consider identifying online basis functions
for a set of $\mu_j$'s following the analysis in \cite{chung2015online}.
This will give a
local eigenvalue problem. Another important problem for our future
consideration is  to identify the values of
$\mu_1$,..., $\mu_J$ by adaptivity.

\bibliographystyle{plain}
\bibliography{references1}

\end{document}